\documentclass[preprint]{imsart}

\RequirePackage{amsthm,amsmath,amsfonts,amssymb}
\RequirePackage[numbers]{natbib}
\RequirePackage[colorlinks,citecolor=blue,urlcolor=blue]{hyperref}
\RequirePackage{graphicx,csquotes,bm}

\startlocaldefs

\newtheorem{theorem}{Theorem}[section]
\newtheorem{lemma}[theorem]{Lemma}
\newtheorem{proposition}[theorem]{Proposition}
\newtheorem{cor}[theorem]{Corollary}
\newtheorem{rem}{Remark}
\newtheorem{definition}[theorem]{Definition}

\newcommand{\R}{\mathbb{R}}
\newcommand{\E}{\mathbb{E}}
\renewcommand{\P}{\mathbb{P}}

\endlocaldefs

\begin{document}

\begin{frontmatter}
\title{Inference on a class of exponential families on permutations }
\runtitle{Inference on a class of exponential families on permutations}

\begin{aug}
\author[A]{\fnms{Sumit} \snm{Mukherjee}
\ead[label=e1]{sm3949@columbia.edu}},
\and
\author[B]{\fnms{Daiki} \snm{Tagami}\ead[label=e3]{daiki.tagami@hertford.ox.ac.uk}}
\address[A]{Department of Statistics,
Columbia University, New York, USA; 
\printead{e1}}

\address[B]{Department of Statistics,
University of Oxford, Oxford, UK; 
\printead{e3}}
\end{aug}

\begin{abstract}
In this paper we study a class of exponential family on permutations, which includes some of the commonly studied Mallows models. We show that the pseudo-likelihood estimator for the natural parameter in the exponential family is asymptotically normal, with an explicit variance. Using this, we are able to construct asymptotically valid confidence intervals. We also show that the MLE for the same problem is consistent everywhere, and asymptotically normal at the origin. In this special case, the asymptotic variance of the cost effective pseudo-likelihood estimator turns out to be the same as the cost prohibitive MLE. To the best of our knowledge, this is the first inference result on permutation models including Mallows models, excluding the very special case of Mallows model with Kendall's Tau.
\end{abstract}


\end{frontmatter}
\section{Introduction}

Ranking data arises naturally in a variety of applications when one compares between several items.  One of the most well studied class of models on the space of rankings/permutations are the celebrated Mallows models, first introduced in \cite{Mallows}. Since then, Mallows models (and their variants) have received significant attention in statistics (\cite{Cr,CFV,FCo, FV1,FV2,Mar, Mukherjee}), probability (\cite{BB,BDMW,DR,GP,HE21,EJP,SS}), and machine learning (\cite{ABSV,CBBK,KHJ, LL,LM,MB2, MPPB}).
In \cite{Mukherjee}, the authors introduce a class of exponential family models on the space of permutations, which includes some of the commonly studied Mallows models. In these models, estimating the parameters via the maximum likelihood estimator (MLE) is computationally infeasible, as the normalizing constant of the exponential family is very hard to analyze, and current technology only allows for a crude (leading order) asymptotics of the log normalizing constant. To bypass this, they study estimation of the natural parameter in the exponential family, using a pseudo-likelihood based approach (\cite{B1,B2}), and derive rates of consistency of their estimator. However, the question of asymptotic distribution of the pseudo-likelihood estimator (PLE) has remained open. And it is of interest to get our hands on the asymptotic distribution, as that will allow us to do more inferential tasks, such as carrying out tests of hypothesis, and constructing confidence intervals.
\\

An interesting observation is that the  proposed model in \cite[Eq (1.1)]{Mukherjee} (see also equation \eqref{eq:model1} of the current paper) has connections to entropic optimal transport (EOT). In particular, the asymptotics of the log normalizing constant of this model (see \cite[Thm 1.5 (a)]{Mukherjee}) gives rise to the same optimization problem studied in entropic optimal transport, with the added restriction of uniform marginals in the permutation case (see \cite[Eq 5]{HLP} or \cite[Eq 1.1]{GNB}). And this is surprising, as there seems to be no direct connection between the two problems. 
It is possible that tools from EOT will be useful here, but this is beyond the scope of the current paper.
\\

In this paper we propose a multivariate generalization of the exponential family proposed in \cite{Mukherjee}. In this general setting, we study the problem of estimating the multivariate natural parameter, having observed one random permutation from this model. We derive the asymptotic limiting distribution of the multivariate PLE, for all values of the true parameter. In particular, this result applies to the one parameter model studied in \cite{Mukherjee}. We also show that the MLE is consistent at all parameter values. Focusing at the origin, we compute the asymptotic distribution of the MLE, which shows that the asymptotic variance of both the estimators is the same. Thus, at least at the origin, the computationally tractable PLE performs equally well (in terms of asymptotic performance), when compared to the computationally inefficient MLE. We also use the PLE to construct asymptotically valid confidence intervals.
\\

One question which is currently out of scope of the current draft is the asymptotic distribution of the MLE, for all values of the true parameter regime. Even though the MLE is not computable, comparing the errors committed by the MLE and the PLE seems to be of interest. Another direction of potential future research is to study the proposed model of this paper from a semi-parametric point of view, where the form of the underlying sufficient statistic is not known a-priori. In this case, we need to have access to multiple independent permutation samples from the underlying model, as it is clear that having one sample from the model (as in the current draft) will not suffice in a semi-parametric setup.
\\

We will now introduce some notation, which will allow us to state our main results.

\subsection{Notation}

For any positive integer $n$, let $[n]:=\{1,2,\ldots,n\}$. Let $S_n$ denote the set of all permutations of $[n]$.  Given any permutation $\pi\in S_n$, we can write $\pi=(\pi(1),\ldots,\pi(n))$ where $\pi(i)\in [n]$ is the image of $i$ under the permutation $\pi$. Define an exponential family on $S_n$ via the following p.m.f.
\begin{align}\label{eq:model}
\P_{n,\bm{\theta}}(\pi)=\exp\Big({\bm \theta}^{\top}{\bf T}(\pi)-Z_n(\bm{\theta})\Big).
\end{align}
Here 
\begin{enumerate}
\item[(i)]
${\bf T}(\pi):=\sum_{i=1}^n{\bf f}\Big(\frac{i}{n},\frac{\pi(i)}{n}\Big)$ is the sufficient statistic, where
${\bf f}=(f_1,\ldots,f_L):[0,1]^2\mapsto \R^L$ is a vector valued continuous function;

\item[(ii)]
$\bm{\theta}\in \R^L$ is a vector valued real parameter;

\item[(iii)]
$Z_n(\bm{\theta})$ is the log normalizing constant, given by
\begin{align}\label{eq:log_norm}
Z_n(\bm{\theta}):=\log\left[ \sum_{\pi \in S_n} \exp\Big( \sum_{i=1}^n{\bm \theta}^{\top} {\bf f}\Big(\frac{i}{n},\frac{\pi(i)}{n}\Big)\Big)\right].
\end{align}
\end{enumerate}
Throughout this paper we will assume that ${\bf f}$ is completely specified, and we will focus on inference about the unknown parameter $\bm{\theta}$. Since ${\bf f}$ is fixed throughout this paper, we have removed the dependence of ${\bf f}$ in \eqref{eq:model}, and we follow this convention throughout the rest of this paper for simplicity of notation.
It is impossible to estimate ${\bm \theta}$, if the functions $(f_1,\ldots,f_L)$ are linearly dependent, as in that case the family $\P_{n,{\bm \theta}}$ is not identifiable in ${\bm \theta}$. Thus throughout we will assume that
${\bf f}$ is linearly independent, i.e.~there does not exist a vector ${\bf d}\in \R^L$ such that ${\bf d}^{\top}{\bf f}\stackrel{a.s.}{=}0$.
\\

Another important observation is that if we replace the function $f_r(x,y)$ by $\tilde{f}_r(x,y)=f_r(x,y)+a_r(x)+b_r(y)$ for arbitrary functions $a_r(.), b_r(.)$, then the underlying model remains unchanged, as
$$\sum_{i=1}^n\tilde{f}_r\Big(\frac{i}{n},\frac{\pi(i)}{n}\Big)=\sum_{i=1}^nf_r\Big(\frac{i}{n},\frac{\pi(i)}{n}\Big)+\sum_{i=1}^n\Big[a_r\Big(\frac{i}{n}\Big)+b_r\Big(\frac{i}{n}\Big)\Big].$$
This allows us to restrict our function class via the following assumption.
\begin{definition}\label{assB}
Let $\mathcal{C}$ denote the collection of all continuous functions $F:[0,1]^2\mapsto \R$ such that
\begin{align}\label{eq:row_col}
\int_0^1F(.,z)dz=0\text{ and } \int_0^1F(z,.)dz=0,
\end{align}
and $F$ is not $0$ a.s.~If ${\bf f}=(f_1,\ldots,f_L)$ is a vector function, we will say ${\bf f}\in \mathcal{C}$, if $f_r\in \mathcal{C}$ for all $r\in [L]$.
\end{definition}
We are now ready to state our main results.

\subsection{Main results}

As indicated above, the log normalizing constant $Z_n(\bm{\theta})$ in \eqref{eq:log_norm} is not available in closed form, and approximating the normalizing constant using numerical techniques/MCMC based methods is a well known challenging problem. 
Using \cite[Thm 1.5]{Mukherjee}) we will now give an asymptotic estimate for $Z_n(\bm{\theta})$. 
\begin{definition}\label{def:sumit}

Let $\mathcal{M}$ denote the space of all probability measures on $[0,1]^2$ with uniform marginals, and let $u\in \mathcal{M}$ be the uniform distribution on $[0,1]^2$. For any measure $\mu$ on the unit square and a bounded measurable function $F:[0,1]\mapsto \R$, let $\mu(F):=\int_{[0,1]^2}Fd\mu$. Finally, let $D(.||.)$ denote the Kullback-Leibler divergence between two probability measures on the same space (typically $[0,1]^2$ or $S_n$).

\end{definition}

\begin{proposition}\label{prop:sumit}
\begin{enumerate}
\item[(a)]
Suppose $Z_n({\bm \theta})$ is as defined in \eqref{eq:log_norm}, and ${\bf f}$ is continuous. Then we have
$$\lim_{n\to\infty}\frac{Z_n({\bm \theta})-\log n!}{n}=\sup_{\mu\in \mathcal{M}}\Big\{{\bm \theta}^{\top} \mu({\bf f}) -D(\mu||u)\Big\}=:Z({\bm \theta}).$$

\item[(b)]
The supremum of part (a) is attained at a unique measure $\mu_{{\bm \theta}}\in \mathcal{M}$ which has a strictly positive density $\rho_{{\bm \theta}}$ (say) with respect to Lebesgue measure. 

\item[(c)]
Under $\P_{n,{\bm \theta}}$ the sequence of empirical measures $\frac{1}{n}\sum_{i=1}^n\delta_{\frac{i}{n},\frac{\pi(i)}{n}}$ converges weakly in probability to the measure $\mu_{{\bm \theta}}$. Consequently, 
$$\frac{{\bf T}(\pi)}{n}\stackrel{P}{\to} \mu_{{\bm \theta}}({\bf f})=:{\bf z}({\bm \theta}).$$

\item[(d)]
The  function ${\bf z}(.):\R^L\mapsto \R^L$ defined in part (c) is continuous on $\R^L$.

\item[(e)]
The function $Z(.)$ defined in part (a) is differentiable  in $\R^L$, with $$\nabla Z({\bm \theta})={\bf z}({\bm \theta})=\mu_{{\bm \theta}}({\bf f}).$$

\item[(f)]
If ${\bf f}$ is linearly independent, and ${\bf f}\in \mathcal{C}$, (see Definition \ref{assB}), then the function $Z(.)$ is strictly convex, i.e.~for any ${\bm \theta}_0\ne {\bm \theta}_1$, we have
\begin{align}\label{eq:strict_convex}
Z({\bm \theta}_1)>Z({\bm \theta}_0)+({\bm \theta}_1-{\bm \theta}_0)^{\top}{\bf z}({\bm \theta}_0).
\end{align}
\end{enumerate}
\end{proposition}

\begin{rem}
It is worthwhile to note that strict convexity of the limiting log normalizing constant shown in Proposition \ref{prop:sumit} part (f) above directly translates into asymptotic identifiability of the model. This will be utilized during the proof of Theorem \ref{thm:ml} part (a), where we show consistency of the MLE. This also follows on using \eqref{eq:kl} to observe that
$$Z({\bm \theta}_1)-Z({\bm \theta}_0)-({\bm \theta}_1-{\bm \theta}_0)^{\top }{\bf z}( {\bm \theta_0})=\lim_{n\to\infty}\frac{1}{n}D(\P_{n,{\bm \theta_0}}||\P_{n,{\bm \theta}_1}),$$
where $D(.||.)$ is the Kullback-Leibler divergence as in Definition \ref{def:sumit}.
\end{rem}
We now introduce the pseudo-likelihood estimator (PLE) for ${\bm \theta}$. 
 
\begin{definition}\label{def:pl}

Let $\mathcal{E}_n:=\{(i,j):1\le i<j\le n\}$. For any $(i,j)\in \mathcal{E}_n$, the conditional distribution of $\pi(i)$ and $\pi(j)$ given $\{\pi(k),k\ne i,j\}$ is given by
\begin{align}\label{eq:pi_condition}
\notag&\P_{n,\bm{\theta}}\Big(\pi(i)=\sigma(i), \pi(j)=\sigma(j)|\pi(k)=\sigma(k), k\ne i,j\Big)\\
=&\frac{\exp\Big({\bm \theta}^{\top} {\bf f}\Big(\frac{i}{n}, \frac{\sigma(i)}{n}\Big)+{\bm \theta}^{\top}{\bf f}\Big(\frac{j}{n}, \frac{\sigma(j)}{n}\Big)\Big)}{\exp\Big({\bm \theta}^{\top} {\bf f}\Big(\frac{i}{n}, \frac{\sigma(i)}{n}\Big)+{\bm \theta}^{\top}{\bf f}\Big(\frac{j}{n}, \frac{\sigma(j)}{n}\Big)\Big)+\exp\Big({\bm \theta}^{\top} {\bf f}\Big(\frac{i}{n}, \frac{\sigma(j)}{n}\Big)+{\bm \theta}^{\top}{\bf f}\Big(\frac{j}{n}, \frac{\sigma(i)}{n}\Big)\Big)}\nonumber\\
&=\frac{\exp\Big({\bm \theta}^{\top} {\bf y}_{\sigma}(i,j)\Big)}{1+\exp\Big({\bm \theta}^{\top} {\bf y}_{\sigma}(i,j)\Big)},
\end{align}
where $${\bf y}_{\sigma}(i,j):={\bf f}\Big(\frac{i}{n},\frac{\sigma(i)}{n}\Big)+{\bf f}\Big(\frac{j}{n},\frac{\sigma(j)}{n}\Big)-{\bf f}\Big(\frac{i}{n},\frac{\sigma(j)}{n}\Big)-{\bf f}\Big(\frac{j}{n},\frac{\sigma(i)}{n}\Big)$$ for $\sigma\in S_n$. 
Then the pseudo-likelihood of $\bm{\theta}$ is defined by the product of these conditional distributions over all pairs $(i,j)\in \mathcal{E}_n$, i.e.
$$PL_n(\pi,\bm{\theta})=\prod_{(i,j)\in \mathcal{E}_n} \frac{\exp\Big({\bm \theta}^{\top} {\bf y}_{\pi}(i,j)\Big)}{1+\exp\Big({\bm \theta}^{\top} {\bf y}_{\pi}(i,j)\Big)},$$
Then the pseudo-likelihood estimator (PLE) is obtained by solving the equation
\begin{align*}
{\bf L}_n(\pi,\bm{\theta})=\nabla \log PL_n(\pi,\bm{\theta})=0.
\end{align*}
A direct calculation gives
\begin{align}\label{eq:pldef}
{\bf L}_n(\pi,\bm{\theta})=\nabla \log PL_n(\pi,\bm{\theta})=\sum_{(i,j)\in \mathcal{E}_n}\frac{1}{1+e^{{\bm \theta}^{\top} {\bf y}_{\pi}(i,j)}} {\bf y}_{\pi}(i,j).
\end{align}
\end{definition}
It is shown in \cite[Thm 1.11]{Mukherjee} that if $L=1$ (i.e.~${\bm \theta}=\theta\in \R$), the expression \eqref{eq:pldef} has a unique zero $\hat{\theta}_{PL}$ with probability tending to $1$, which satisfies $\sqrt{n}(\hat{\theta}_{PL}-\theta)=O_P(1)$.  The question of asymptotic distribution of $\hat{\bm \theta}_{PL}$ has remained open, even in the one dimensional case. Our first main result addresses this, by solving the more general multivariate analogue of this problem. 
%
%
%

\begin{theorem}\label{thm:pl}
Suppose $\pi$ is a random permutation from the model \eqref{eq:model}, where ${\bf f}$ is linearly independent, and ${\bf f}\in \mathcal{C}$.  Then, denoting the true parameter as $\bm{\theta}_0\in \R^L$, the following conclusions hold:
\begin{enumerate}
\item[(a)]
With ${\bf L}_n(\pi,\bm{\theta})$ as defined in \eqref{eq:pldef}, we have
$$n^{-3/2}{\bf L}_n(\pi,\bm{\theta}_0)\stackrel{D}{\to} N\Big(0,\Sigma({\bm \theta}_0)\Big),$$
where 
$$\Sigma_{pq}({\bm \theta}):= \int_{[0,1]^6}\frac{g_p({\bf z}_1,{\bf z}_2)g_q({\bf z}_1,{\bf z}_3)}{(1+e^{{\bm \theta}^{\top}{ \bf g}({\bf z}_1,{\bf z}_2)})(1+e^{{\bm \theta}^{\top}{ \bf g}({\bf z}_1,{\bf z}_3)})}  \prod_{a=1}^3 d\mu_{{\bm \theta}}({\bf z}_a) $$
for $p,q\in [L]$. Here the measure $\mu_{{\bm \theta}}$ is as in Proposition \ref{prop:sumit}, and ${\bf g}:[0,1]^4\mapsto \R^L$ is defined by 
\begin{align*}
{\bf g}((x_1,y_1),(x_2,y_2)):={\bf f}(x_1,y_1)+{\bf f}(x_2,y_2)-{\bf f}(x_1,y_2)-{\bf f}(x_2,y_1).
\end{align*}

\item[(b)]
The vector equation ${\bf L}_n(\pi,\bm{\theta})=\bm{0}$ has a unique root $\bm{\hat{\theta}_{PL}}$ with probability tending to $1$, which satisfies $$\sqrt{n}(\hat{\bm \theta}_{PL}-{\bm \theta}_0)\stackrel{D}{\to}N\Big(0,[A({\bm \theta}_0)]^{-1}\Sigma({\bm \theta}_0) [A({\bm \theta}_0)]^{-1}\Big).$$
Here $A({\bm \theta})$ is a positive definite $L\times L$ matrix defined by $$A_{pq}({\bm \theta}):=\frac{1}{2}\int_{[0,1]^4} \frac{g_p({\bf z}_1,{\bf z}_2)g_q({\bf z}_1,{\bf z}_2) e^{{\bm \theta}^{\top} {\bf g}({\bf z}_1,{\bf z}_2)}}{(1+e^{{\bm \theta}^{\top}{\bf  g}({\bf z}_1,{\bf z}_2)})^2} \prod_{a=1}^2 d\mu_{{\bm \theta}}({\bf z}_a) .$$

\end{enumerate}
\end{theorem}

As a special case of Theorem \ref{thm:pl}, setting $L=1$ we immediately get the following corollary.

\begin{cor}\label{cor:pl}
Suppose $\pi$ is a random permutation from the model 
\begin{align}\label{eq:model1}
\P_{n,\theta}(\pi)=\exp\Big(\theta \sum_{i=1}^nF\Big(\frac{i}{n}, \frac{\pi(i)}{n}\Big)-Z_n(\theta)\Big),
\end{align} where $F\in \mathcal{C}$ and $\theta\in \R$. Then, denoting the true parameter as ${\theta}_0\in \R$, the PLE $\hat{\theta}_{PL}$ exists with probability tending to $1$, and satisfies
%
$$\sqrt{n}(\hat{ \theta}_{PL}-{ \theta}_0)\stackrel{D}{\to}N\Big(0,\frac{\sigma^2(\theta_0)}{a(\theta_0)}\Big),$$
where 
\begin{align*}
\sigma^2(\theta):=& \int_{[0,1]^6}\frac{G({\bf z}_1,{\bf z}_2)G({\bf z}_1,{\bf z}_3)}{(1+e^{\theta G({\bf z}_1,{\bf z}_2)})(1+e^{\theta G({\bf z}_1,{\bf z}_3)})}  \prod_{a=1}^3 d\mu_{{ \theta}}({\bf z}_a) ,\\
a(\theta):=&\frac{1}{2}\int_{[0,1]^4} \frac{G^2({\bf z}_1,{\bf z}_2) e^{\theta G({\bf z}_1,{\bf z}_2)}}{(1+e^{\theta G({\bf z}_1,{\bf z}_2)})^2} \prod_{a=1}^2 d\mu_{{ \theta}}({\bf z}_a).
\end{align*}
Here $$G((x_1,y_1),(x_2,y_2))=F(x_1,y_1)+F(x_2,y_2)-F(x_1,y_2)-F(x_2,y_1),$$
and $\mu_{\theta}$ is as in Proposition \ref{prop:sumit} part (c) with $L=1$.


\end{cor}
\begin{rem}
We point out here that taking $F(x,y)=-|x-y|$ in \eqref{eq:model1} we get the Mallows model with Spearman's Footrule as sufficient statistic, and taking $F(x,y)=-(x-y)^2$ (or $F(x,y)=xy$) in \eqref{eq:model1} we get the Mallows model with Spearman's Rank Correlation as sufficient statistic.  We refer to \cite{diaconis,Mallows,Mukherjee} for more background on these Mallows models, as well as other Mallows models considered in the literature.
To the best of our knowledge, Theorem \ref{thm:pl} and Corollary \ref{cor:pl} are the first results which can allow inferential procedures such as testing of hypothesis in the model \eqref{eq:model}, which includes these Mallows models. Prior to our work, the only inference results that we are aware of is in the Mallows model with Kendall's Tau as sufficient statistic, where the model has an explicit normalizing constant, and hence is very tractable. 
\end{rem}

A natural question is whether the PLE is asymptotically optimal, in the sense that it has the smallest asymptotic variance. Even though the MLE is incomputable, it is still expected to be the \enquote{gold standard} in terms of estimators for statistical efficiency, at least for nice exponential families such as \eqref{eq:model}. Thus one may ask whether one can compare the performance of the MLE to that of the PLE. Towards this direction, our next result shows that the MLE is consistent for all ${\bm \theta}_0\in \R^L$. Analyzing the asymptotic distribution of the MLE is more delicate, and requires precise asymptotic properties of the log normalizing constant.  We are able to carry out this program in a neighborhood of the origin, thus capturing the CLT of the MLE at $\bm{\theta}_0=\bm{0}$. 

\begin{theorem}\label{thm:ml}
Suppose $\pi$ is a random permutation from the model \eqref{eq:model}, where ${\bf f}=(f_1,\ldots,f_L)$ are linearly independent, ${\bf f}\in \mathcal{C}$, and the true parameter is $\bm{\theta}_0\in \R^L$. Then the following conclusions hold:

\begin{enumerate}
\item[(a)]
The equation ${\bf T}(\pi)=\nabla Z_n(\bm{\theta})$ has a unique root $\bm{\hat{\theta}_{ML}}$ with probability tending to $1$, which satisfies
$$\hat{\bm \theta}_{ML}\stackrel{P}{\rightarrow}{\bm \theta}_0.$$

\item[(b)]
If ${\bm \theta}_0={\bf 0}$, then we have
$$\frac{{\bf T}(\pi)-\nabla Z_n(\bm{0})}{\sqrt{n}}\stackrel{D}{\to}N\Big(\bm{0},\Gamma\Big),\text{ where }
\Gamma_{p,q}:=\int_{[0,1]^2} f_p(x,y)f_q(x,y)dx dy.$$
%

\item[(c)]
If ${\bm \theta}_0={\bf 0}$, then we have
$$\sqrt{n}\bm{\hat{\theta}_{ML}}\stackrel{D}{\to}N\Big(\bm{0},\Gamma^{-1}\Big).$$

\end{enumerate}
\end{theorem}

Comparing Theorems \ref{thm:pl} and \ref{thm:ml}, it follows that the asymptotic distribution of PLE and MLE are both same if ${\bm \theta}_0={\bf 0}$, as shown in the following corollary.

\begin{cor}

With $A( {\bm \theta}), \Sigma({\bm \theta})$ and $\Gamma$ as defined in Theorems \ref{thm:pl} and Theorem \ref{thm:ml}, we have $$A({\bf 0})^{-1} \Sigma({\bf 0})A({\bf 0})^{-1}=\Gamma^{-1}.$$
Consequently, both $\hat{\bm \theta}_{PL}$ and $\hat{\bm \theta}_{ML}$ have the same asymptotic distribution if the true parameter ${\bm \theta}_0={\bm 0}$.

\end{cor}

\begin{rem}
It remains to be seen whether the performance of the PLE matches that of the MLE for all ${\bm \theta}_0\in \R^L$.  The main theoretical bottleneck is an absence of CLT for ${\bf T}(\pi)$ under the exponential family $\P_{n,{\bm \theta}}$, which is an interesting question in its own right, more so due to its connection to Entropic Optimal Transport, as indicated in the introduction. In particular, \cite{HLP} deduces a CLT for a related exponential family on $[0,1]^{2n}$ arising in EOT, and it is possible that the techniques used there may be extended to cover the CLT for ${\bf T}(\pi)$ under the model \eqref{eq:model1}.
\end{rem}
A related inferential question is the construction of confidence intervals. Our next result constructs an asymptotically valid confidence interval for the parameter ${\bf d}^{\top}{\bm \theta}$, where ${\bf d}\ne {\bf 0}$ is a known vector in $\R^L$.
\begin{lemma}\label{lem:CI}
For $1\le i<j\le n$ and $p,q\in [L]$, set
\begin{align*}
\widehat{\Sigma}_{p,q}(\bm \theta)=&\frac{1}{n^3}\sum_{(i,j),(k,\ell) \in \mathcal{E}_n: |(i,j)\cap (k,\ell)|=1}\frac{y_{\pi,p}(i,j)y_{\pi,q}(i,k)}{(1+e^{{\bm\theta}^\top {\bf y}_\pi(i,j)})(1+e^{{\bm\theta}^\top {\bf y}_\pi(i,k)})},\\
\widehat{A}_{p,q}(\bm \theta)=&\frac{1}{n^2}\sum_{{(i,j)\in \mathcal{E}_n}}\frac{y_{\pi,p}(i,j)y_{\pi,q}(i,j)e^{{\bm \theta}^\top {\bf  y}_{\pi}(i,j)}}{(1+e^{{\bm\theta}^\top {\bf y}_\pi(i,j)})^2},
\end{align*}
where ${\bf y}_\pi(i,j)$ is as in Definition \ref{def:pl}.
Then for any $\alpha\in (0,1)$ and ${\bf d}\in \R^L$ we have
\begin{align*}
\lim_{n\to\infty}\P_{n,{\bm \theta}_0}\left(\Big|{\bf d}^{\top}\hat{\bm \theta}_{PL}-{\bf d}^{\top}{\bm \theta}_0\Big|\le \frac{z_{\frac{\alpha}{2}}}{\sqrt{n}}\cdot \sqrt{{\bf d}^{\top}\Big[[\widehat{A}(\hat{\bm \theta}_{PL})]^{-1}\widehat{\Sigma}(\hat{\bm \theta}_{PL})[\widehat{A}(\hat{\bm \theta}_{PL})]^{-1}\Big]{\bf d}}\right)= 1-\alpha,
\end{align*}
where $z_\alpha$ is the $(1-\alpha)^{th}$ quantile of the standard normal distribution.

\end{lemma}

\subsection{Simulation results}
To illustrate our results, we focus on a specific example of model \eqref{eq:model1} with $F(x,y)=xy$. This 
corresponds to the Mallows model with Spearman's rank correlation. In this case the limiting measure $\mu_{\theta}$ has a density with respect to Lebesgue measure on $[0,1]^2$, of the form
$$e^{\theta xy+a_\theta(x)+a_\theta(y)},$$
where the extra symmetry is due to the fact that $F$ is symmetric in $(x,y)$ (see \cite[Sec 2]{Mukherjee}).
The function $a_\theta(.)$ is uniquely determined by the requirement that the above density has uniform marginals. 
\\

 To illustrate the results of the current paper, we need to be able to simulate from this model efficiently. 
  In \cite{AD} the authors derive an auxiliary variable/hit and run algorithm to simulate from this model, 
which is explained below:
\begin{enumerate}
\item[(i)] Simulate permutation $\pi$ uniformly at random from $S_n$.
\item[(ii)] Given $\pi$, simulate mutually independent random variables $\{U_i\}_{i=1}^n$, with $U_i$ uniform on $[0,e^{(\theta/n^2)i\pi(i)}]$.
\item[(iii)] Given $\{U_i\}_{i=1}^n$, set $b_j:=\max\{(n^2/\theta j)\log{U_j},1\}$. Choose an index $i_1$ uniformly at random from the set $\{j\in [n]:b_j\leq 1\}$, and set $\pi(i_1)=1$. Remove this index from $[n]$ and choose an index $i_2$ uniformly at random from the set $\{j\in [n]:b_j\leq 2\}-\{i_1\}$, and set $\pi(i_2)=2$. In general, having defined indices $\{i_1,\dots,i_{l-1}\}$, remove them from $[n]$. Choose $i_\ell$ uniformly at random from the set of indices $\{j\in [n]:b_j\leq \ell\}-\{i_1,\dots,i_{l-1}\}$ and set $\pi(i_\ell)=\ell$.
\item[(iv)] Iterate between the steps 2 and 3 until convergence.
\end{enumerate}
For our simulations, we iterate steps 2 and 3 a total of 10 times to obtain a single permutation $\pi$. 
Given the permutation $\pi$, we compute the PLE using the bisection method. 
To examine the distribution of the PLE, we repeat the above process $2000$ times, with permutation sizes $n=500,2000,8000$ and true parameter $\theta=2$. The three histograms obtained from this are superimposed in Figure \ref{fig:hist} below, in colors pink, green and blue respectively.
\begin{figure}[h]
\centering
\centering
\includegraphics[width=4in,height=2.8in]
    {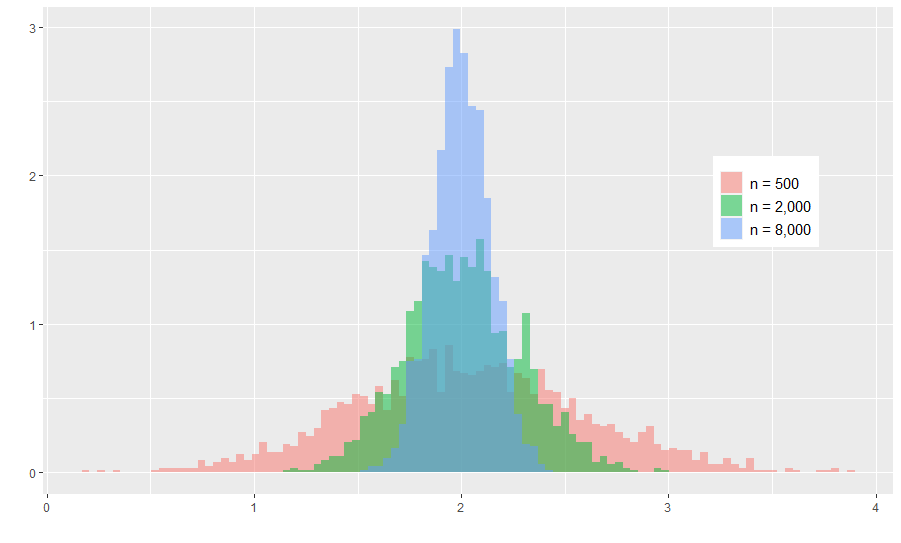}\\
\caption{Histogram of the PLE  in the Spearman's rank correlation model, with permutation size $n=500, 2000, 8000$ and $\theta_0=2$, based on $2000$ samples. }
\label{fig:hist}\vspace{-0.15in}
\end{figure}
From the figure, we see that all the three histograms are bell shaped with a center around $2$. This agrees with Theorem \ref{thm:pl}, which shows that the PLE is asymptotically normal. Also, the width of the histogram decreases as $n$ increases, which agree with the fact that the variance of the PLE is decreasing with $n$.
\\

To compare the performance of the PLE with that of the MLE, we focus on the same model, but set $\theta_0=0$ instead. Note that computing the MLE is a challenge, so we use the following heuristic approximation:



Using \eqref{eq:logpart2} in the proof of part (c) of Theorem \ref{thm:ml}, it follows that if ${ \theta}_0={ 0}$, for any $u\in \R$ we have
$$Z_n'\left(\frac{u}{\sqrt{n}}\right)\approx Z_n'(0)+u\gamma \sqrt{n}, \quad \gamma=\int_{[0,1]^2}F^2(x,y)dx dy.$$
Setting $u=U_n=\sqrt{n}{ \theta}_{ML}$ using $\sqrt{n}$ consistency of MLE under ${ \theta}_0={ 0}$ we have that $U_n=O_P(1)$, and so the above display gives
$${ T}(\pi)=Z_n'\Big(\hat{ \theta}_{ML}\Big)=Z_n'\Big(\frac{U_n}{\sqrt{n}}\Big)\approx Z_n'(0)+n\gamma \hat{ \theta}_{ML},$$
where $T(\pi)=\sum_{i=1}^nF\Big(\frac{i}{n},\frac{\pi(i)}{n}\Big).$
This in turn gives the simple approximation 
\begin{align}\label{eq:mle_approx}
\hat{ \theta}_{ML}\approx \frac{{ T}(\pi)-Z_n'(0)}{n\gamma},
\end{align}
which is very easy to compute. In our example we have 
\begin{align*}
Z_n'(0)=\E_{n,{ \theta}={ 0}}{T}(\pi)=\sum_{i=1}^n\frac{i}{n^2}\E_{n,{ \theta}={ 0}}[\pi(i)]=\frac{(n+1)^2}{4n},\qquad
\gamma=\int_{[0,1]^2} x^2 y^2 dx dy=\frac{1}{9}.
\end{align*}
Plugging these values, one can compute an approximation to the MLE using \eqref{eq:mle_approx}. We note that technically this is not an estimator, as we used the knowledge of $\theta_0=0$ to justify the heuristic.
\\

\begin{figure}[h]
\centering
\centering
\includegraphics[width=4in,height=2.8in]
    {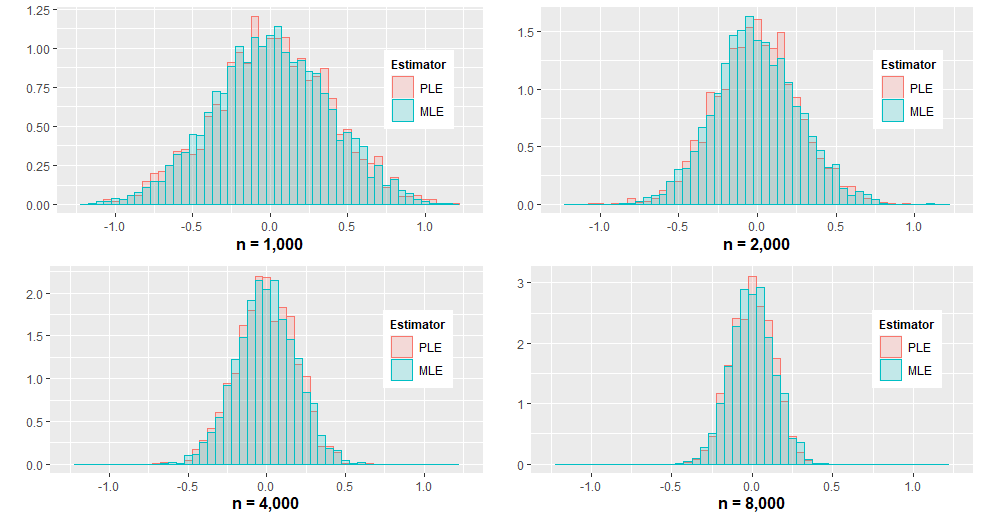}\\
\caption{Histogram of MLE and PLE  in the Spearman's rank correlation model, with permutation size $n=1000, 2000, 4000, 8000$ and $\theta_0=0$, based on $2000$ samples. }
\label{fig:hist2}
\vspace{-0.15in}
\end{figure}

%
%
%
%
In Figure \ref{fig:hist2}, we compare the histogram of this approximate MLE and the PLE where the underlying true parameter is $\theta_0=0$ (i.e.~$\pi$ is a uniformly random permutation). We simulated 2000 permutations $\pi$ with sizes $n=1000,2000,4000,8000$, thereby producing four such comparative histograms (one per value of $n$). 
In each of the above figures, we see that both the histograms overlap, thereby demonstrating that the asymptotic distribution of the MLE and the PLE is the same, when the underlying true parameter is set to ${\bm \theta}_0={\bf 0}$. Also, as expected, we see that the histograms shrink in their width, as the size of the permutation $n$ increases. 
\\

\begin{figure}[h]
\centering
\centering
\includegraphics[width=4in,height=2.2in]
    {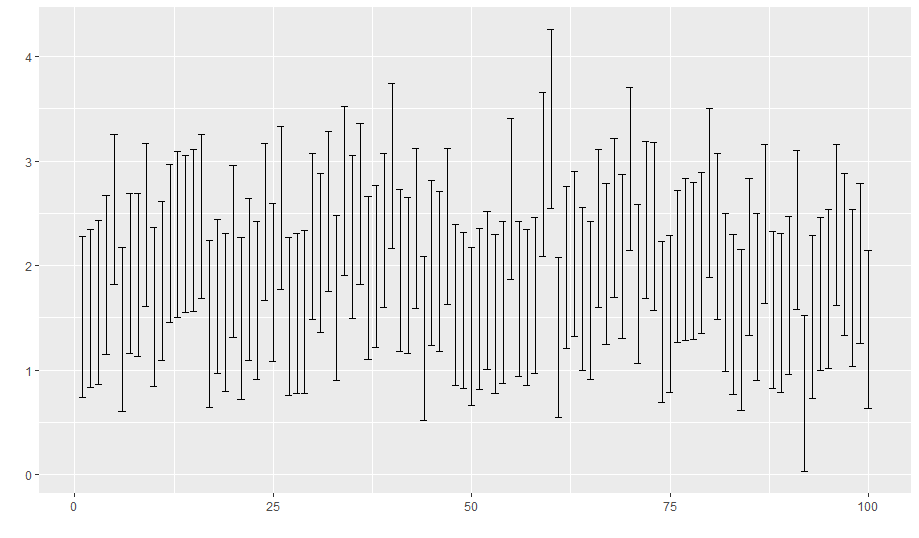}\\
\caption{Confidence Interval based on the PLE  in the Spearman's rank correlation model, with permutation size $n=1000$ and $\theta_0=2$, based on $100$ separate samples.  }
\label{fig:ci}
\vspace{-0.15in}
\end{figure}

Finally, we construct a 95\% confidence interval for the parameter $\theta_0$, when the underlying true parameter is $\theta_0=2$, and the size of the permutation is $n=1000$. The construction of each confidence interval is based on one permutation, so essentially we are doing inference based on one sample (which is of size $n$, of course)! We repeat the construction of this confidence interval $100$ times, and plot the confidence intervals in Figure \ref{fig:ci}.
As it turns out, out of the $100$ confidence intervals so constructed, $95$ of them contain the true value $\theta_0=2$. This validates Lemma \ref{lem:CI}, which claims that this confidence interval has approximate coverage $95\%$ for $n$ large.


%
%
%
%
%
%
%
%

\subsection{Outline of the paper}

In section \ref{sec:main} we prove the main results of this paper, stated in the introduction. The proof of these results are done using supporting lemmas, whose proofs are deferred to section \ref{sec:support}.

\section{Proof of main results}\label{sec:main}

We begin with following three lemmas, whose proof we defer to section \ref{sec:support}. The first two lemmas compute moments of linear functions of ${\bf L}_n(\pi,{\bm \theta})$. These will be used in the proof of Theorem \ref{thm:pl}.

\begin{lemma}\label{lem:conditional_mean}
For $(i,j)\in \mathcal{E}_n$ and for arbitrary ${\bf d}\in \R^L$, setting $${\bf C}_{ij}(\pi):=\frac{1}{1+\exp\Big({\bm \theta}_{0}^{\top} {\bf y}_{\pi}(i,j)\Big)}{\bf y}_{\pi}(i,j)\in \R^L,\text{ and }C_{ij}(\pi,{\bf d})={\bf d}^{\top}{\bf C}_{ij}(\pi),$$
 we have
$$ \E_{n,{\bm \theta}_0}  \Big(C_{ij}(\pi,{\bf d})\Big| \pi(\ell), \ell\ne i,j\Big) =0.$$
\end{lemma}

\begin{lemma}\label{lem:moment}
Let $\{D_{ij},(i,j)\in \mathcal{E}_n\}$ be jointly centered multivariate Gaussian, with $Cov(D_{ij}, D_{k\ell})=\E_{n,{\bm \theta}_0} C_{ij}(\pi,{\bf d}) C_{k\ell}(\pi,{\bf d})$. Thus, setting $D_n:=\sum_{(i,j)\in \mathcal{E}_n} D_{ij}$, 
for any $k\in \mathbb{N}$ we have
$$\lim_{n\to\infty}\frac{\E_{n,{\bm \theta}_0} \Big(\sum_{(i,j)\in \mathcal{E}_n}C_{ij}(\pi,{\bf d})\Big)^k-\E  D_n^k}{n^{\frac{3k}{2}}}=0.$$

\end{lemma}

The third lemma proves a general asymptotic normality for M estimators (which is also applicable in dependent settings), and will be used to prove asymptotic normality of both the PLE (in Theorem \ref{thm:pl} part (b)) and the MLE (in Theorem \ref{thm:ml} part (c)). Even though the lemma is stated for functions over $S_n$ for easy applicability, the proof applies verbatim when $S_n$ is replaced by a more general sample space.

\begin{lemma}\label{lem:consist}
%
%
%
%
%
%

Suppose $\mathcal{S}:S_n\times \R^L\mapsto \R$ be a function, such that for every $\pi\in S_n$ the map $\mathcal{S}(\pi,.)$ is  $\mathcal{C}^2$, and the Hessian $H(\pi,{\bm \theta})\in \R^{L\times L}$ is non-negative definite. 
Fixing ${\bm \theta}_0\in \R^L$, assume further that the following conditions hold under $\P_{n,{\bm \theta}_0}$:

(i) The gradient $\nabla \mathcal{S}(\pi, {\bm \theta})$ satisfies \[n^{-1/2}\nabla \mathcal{S}(\pi, {\bm \theta}_0)\stackrel{D}{\to}N({\bf 0},B_1).\]

(ii) For any sequence ${\bf h}_n\in \R^L$ converging to ${\bf h}$ the Hessian $H(\pi,.)$ satisfies
\[\frac{1}{n}H(\pi,{\bm \theta}_0+n^{-1/2}{\bf h}_n)\to B_2,\]
with $B_2$ positive definite, free of ${\bf h}$. Then, with probability tending to $1$ there exists a unique $\hat{\bm \theta}$ which satisfies $\nabla \mathcal{S}(\pi, {\bm \theta})={\bf 0}$. Further we have
\[ \sqrt{n}(\hat{\bm\theta}-{\bm \theta}_0)\stackrel{D}{\rightarrow}N({\bf 0}, B_2^{-1}B_1B_2^{-1}).\]

\end{lemma}

In a similar spirit, the fourth lemma gives a general consistency result for M estimators, which is also applicable in dependent settings.  We will use it to prove consistency of the MLE in \ref{thm:ml} part (a). As above,  the lemma can be easily extended to cover settings where the sample space is not $S_n$.

\begin{lemma}\label{lem:consistb}

 Let $\mathcal{S}:S_n\times \R^L\mapsto \R$ be a function which is strictly convex in the second argument.
 Let ${\bm \theta}_0\in \R^L$ be such that the following conditions hold under $\P_{n,{\bm \theta}_0}$:
 
 (i) For every ${\bm \theta}_1\ne {\bm \theta}_0$ there exists $\varepsilon:=\varepsilon_{\theta_1}>0$ such that
 $$\P_{n,{\bm \theta}_0}\Big(\mathcal{S}(\pi, {\bm \theta}_0)+\varepsilon<\mathcal{S}(\pi, {\bm \theta}_1)\Big)=1.$$
 
 (ii) There exists a random variable $M(\pi)$ which is $O_P(1)$, such that for all ${\bm \theta}_1,{\bm \theta}_2\in \R^L$ we have 
$$|\mathcal{S}(\pi,{\bm \theta}_1)-\mathcal{S}(\pi,{\bm \theta}_2)|\le M(\pi)\|{\bm \theta}_1-{\bm \theta}_2\|_2.$$

Then there exists a unique minimizer $\hat{\bm \theta}_n$ of the function $\mathcal{S}(\pi,.)$ in $\R^L$ with probability tending to $1$. Further, we have
$$\hat{\bm \theta}_n\stackrel{P}{\to}{\bm \theta}_0.$$

\end{lemma}

Armed with these lemmas, we now complete  the proofs of Theorem \ref{thm:pl}, Theorem \ref{thm:ml}, and Lemma \ref{lem:CI}.

\subsection{Proof of Theorem \ref{thm:pl}}

\begin{enumerate}
\item[(a)]
For part (a), it suffices to show that for any ${\bf d}\in \R^L$ we have
\begin{align*}
\frac{1}{n\sqrt{n}}{\bf d}^{\top}{\bf L}_{n}(\pi,{\bm \theta}_0)\stackrel{D}{\to} N\left(0,\bm{d}^{\top}\Sigma({\bm \theta}_0)\bm{d}\right).
\end{align*}
With $C_{ij}(\pi,{\bf d})$ as defined in Lemma \ref{lem:conditional_mean}, we have
$${\bf d}^{\top}{\bf L}_{n}(\pi,\bm{\theta}_0)=\sum_{(i,j)\in \mathcal{E}_n}C_{ij}(\pi,{\bf d}).$$
%
%
%
Invoking Lemma \ref{lem:moment}, it suffices to show that $\frac{D_n}{n\sqrt{n}}\stackrel{D}{\to}N(0,\bm{d}^{\top}\Sigma\bm{d})$. Since $D_n$ is a centered Gaussian by construction, it suffices to show that
\begin{align}\label{eq:claim2}
\frac{Var_{n,{\bm \theta}_0}(D_n)}{n^3}=\frac{1}{n^3}\sum_{(i,j),(k,\ell)\in \mathcal{E}_n} \E_{n,{\bm \theta}_0} C_{ij}(\pi,{\bf d})C_{k\ell}(\pi,{\bf d})\to \bm{d}^{\top}\Sigma({\bm \theta}_0)\bm{d}.
\end{align}

Proceeding to verify the display above, 
combining \cite[Cor 1.12]{EJP} and \cite[Thm 1.4]{EJP} we conclude that for any $\ell\in \mathbb{N}$ we have
\begin{align}\label{eq:ejp}
\lim_{n\to\infty} \sup_{{\bf p}, {\bf q}\in S(n,\ell)} \left|\frac{n^\ell \P_{n,{\bm \theta}_0}\Big(\pi(p_a)=q_a,a\in [\ell]\Big)}{\prod_{a=1}^\ell \rho_{{\bm \theta}}\Big(\frac{p_a}{n}, \frac{q_a}{n}\Big)}-1\right|=0.
\end{align}
In the display above, $S(n,\ell)\subset [n]^\ell$ is the set of all $\ell$ tuples which consists of distinct entries, and $\rho_{{\bm \theta}}$ is as defined in Proposition \ref{prop:sumit}. Essentially, \eqref{eq:ejp} implies that for any finite collection of indices $\{p_1,\ldots,p_\ell\}$, the random variables $\{\pi(p_1),\ldots,\pi(p_\ell)\}$ are approximately independent, and identifies their joint limiting distribution. Using Lemma \ref{lem:conditional_mean}, it follows that if $|(i,j)\cap(k,\ell)|=0$, then
\begin{align}\label{eq:c4}
\E_{n,{\bm \theta}_0} C_{ij}(\pi,{\bf d})C_{k\ell}(\pi,{\bf d})=0.
\end{align}
Also, it is immediate that
\begin{align}\label{eq:c3}
\Big|\sum_{(i,j),(k,\ell)\in \mathcal{E}_n:|(i,j)\cap (k,\ell)|= 2}\E_{n,{\bm \theta}_0} C_{ij}(\pi,{\bf d}) C_{k\ell}(\pi,{\bf d})\Big|=O(n^2).
\end{align}
Thus the leading contribution in display \eqref{eq:claim2} comes from terms of the form $\E_{n,{\bm \theta}_0} C_{ij}(\pi,{\bf d}) C_{k\ell}(\pi,{\bf d})$ which satisfy $|(i,j)\cap (k,\ell)|=1$. Invoking \eqref{eq:ejp}, we get
\begin{align}\label{eq:c2}
\notag&\sum_{(i,j),(k,\ell)\in \mathcal{E}_n:|(i,j)\cap (k,\ell)|=1}\E_{n,{\bm \theta}_0} C_{ij}(\pi,{\bf d}) C_{k\ell}(\pi,{\bf d})\\
\to&\sum_{p,q=1}^Ld_pd_q\int_{[0,1]^6}\frac{g_p({\bf z}_1,{\bf z}_2)}{1+e^{{\bm \theta}_0^{\top}{\bf g}({\bf z}_1,{\bf z}_2)}}  \frac{g_q({\bf z}_1,{\bf z}_3)}{1+e^{{\bm \theta}_0^{\top}{\bf g}({\bf z}_1,{\bf z}_3)}}  \prod_{a=1}^3 d\mu_{{\bm \theta}_0}({\bf z}_a)
=\bm{d}^{\top}\Sigma({\bm \theta}_0)\bm{d}.
\end{align}
where ${\bf g}(.,.)$ is as in the statement of Theorem \ref{thm:pl}. Combining \eqref{eq:c4}, \eqref{eq:c3} and \eqref{eq:c2} verifies \eqref{eq:claim2}, and hence completes the proof of part (a).
\\

\item[(b)]


%

For part (b), we will invoke Lemma \ref{lem:consist} with the following choices:
\[ \mathcal{S}(\pi,{\bm \theta})=-n^{-1}\log PL_n(\pi,{\bm \theta}), \quad B_1=\Sigma({\bm \theta}_0), \quad B_2=A({\bm \theta}_0).\]
In this case we have 
\[n^{-1/2}\nabla \mathcal{S}(\pi, {\bm \theta}_0)=-n^{-3/2}{\bf L}_n(\pi,{\bm \theta}_0)\stackrel{D}{\to}N\Big({\bf 0},\Sigma({\bm \theta}_0)\Big)\]
by part (a), and so assumption (i) of Lemma \ref{lem:consist} holds. Taking a second derivative gives
\begin{align*}
H_{pq}(\pi,{\bm \theta})=\frac{\partial^2 \mathcal{S}(\pi,{\bm \theta})}{\partial \theta_p \partial \theta_q}=&n^{-1}\sum_{(i,j)\in \mathcal{E}_n} y_{\pi,p}(i,j)y_{\pi,q}(i,j) \frac{e^{{\bm \theta}^{\top}{\bf  y}_\pi(i,j)}}{(1+e^{{\bm \theta}^{\top}{\bf  y}_\pi(i,j)})^2}.
\end{align*}
Consequently, the Hessian $H(\pi,{\bm \theta})$ is non-negative definite. Further, if $\{{\bf h}_n\}_{n\ge 1}$ is a sequence in $\R^L$ converging to ${\bf h}$, then using \eqref{eq:ejp} we have
\begin{align}\label{eq:CI2}
\notag&\frac{1}{n}H_{pq}(\pi,{\bm \theta}_0+n^{-1/2}{\bf h}_n)\\
\stackrel{P}{\to}&\frac{1}{2}\int_{[0,1]^4} g_p({\bf z}_1,{\bf z}_2)g_q({\bf z}_1,{\bf z}_2) \frac{e^{{\bm \theta}_0^{\top} {\bf g}({\bf z}_1,{\bf z}_2)}}{{(1+e^{{\bm \theta}_0^{\top} {\bf g}({\bf z}_1,{\bf z}_2)})^2}}\ \prod_{a=1}^2 d\mu_{{\bm \theta}_0}({\bf z}_a)=A_{pq}({\bm \theta}_0),
\end{align}
where $A(.)$ is as defined in Theorem \ref{thm:pl}. To verify assumption (ii), it remains to show that $A({\bm \theta})$ is positive definite. 
 To this effect, for any vector ${\bf b},{\bm \theta}\in \R^L$ we have
%
%
%
%
%
\begin{align*}
\nonumber\bm{b}^{\top}A({\bm \theta})\bm{b}&=\frac{1}{2}\sum_{p,q=1}^Lb_pb_q\int_{[0,1]^4} \frac{g_p({\bf z}_1,{\bf z}_2)g_q({\bf z}_1,{\bf z}_2) e^{{\bm \theta}^{\top} {\bf g}({\bf z}_1,{\bf z}_2)}}{(1+e^{{\bm \theta}^{\top} {\bf g}({\bf z}_1,{\bf z}_2)})^2} \prod_{a=1}^2 d\mu_{{\bm \theta}}({\bf z}_a)\\&=\frac{1}{2}
\int_{[0,1]^4}\frac{\Big(\sum_{r=1}^L b_rg_r({\bf z}_1,{\bf z}_2)
\Big)^2e^{{\bm \theta}^{\top} {\bf g}({\bf z}_1,{\bf z}_2)}}{(1+e^{{\bm \theta}^{\top} {\bf g}({\bf z}_1,{\bf z}_2)})^2}\prod_{a=1}^2 d\mu_{{\bm \theta}}({\bf z}_a)
\end{align*}
Clearly $A({\bm \theta})$ is non-negative definite, and to verify positive definiteness it suffices to show the RHS above is not $0$. If this is $0$, then we must have
%
%
\begin{align}\label{eq:as}
\sum_{r=1}^Lb_rg_r({\bf z}_1,{\bf z}_2)\stackrel{a.s.}{=}0
\end{align} with respect to $\mu_{{\bm \theta}^{(0)}}^{\otimes 2}$. Invoking Proposition \ref{prop:sumit} part (b) we have that $\mu_{{\bm \theta}^{(0)}}$ has a strictly positive density with respect to Lebesgue measure on $[0,1]^2$, and so \eqref{eq:as} holds with respect to Lebesgue measure on $[0,1]^4$ as well. Integrating over ${\bf z}_2$ and using the assumption that $f_r\in \mathcal{C}$ gives $$\sum_{r=1}^L b_r f_r({\bf z}_1)\stackrel{a.s.}{=}0$$
with respect to Lebesgue measure on $[0,1]^2$. But this violates the assumption that $(f_1,\ldots,f_r)$ are linearly independent. This contradiction shows that $A({\bm \theta}_0)$ is positive definite, and hence verifies assumption (ii) of Lemma \ref{lem:consist}. The proof of part (b) then follows by invoking the lemma.

\end{enumerate}

\subsection{Proof of Theorem \ref{thm:ml}}

\begin{enumerate}
\item[(a)]

For proving part (a), we will invoke Lemma \ref{lem:consistb} with $$\mathcal{S}(\pi,{\bm \theta})=\frac{1}{n}\Big[Z_n({\bm \theta})-{\bm \theta}^{\top}{\bf T}(\pi)\Big].$$
In this case $\mathcal{S}(\pi,.)$ is convex, as the Hessian equals $\frac{1}{n}Var_{n,{\bm \theta}}({\bf T}(\pi))$. We now claim that $Var_{n,{\bm \theta}}({\bf T}(\pi))$ is not singular, for all large $n$. By way of contradiction, assume that $Var_{n,{\bm \theta}}({\bf T}(\pi))$ is singular for some $n$. Since the two laws $\P_{n,{\bm \theta}}$ and $\P_{n,{\bf 0}}$ are mutually absolutely continuous, it follows that $Var_{n, {\bf 0}}({\bf T}(\pi))$ is singular as well.  But we argue directly in part (b) (see \eqref{eq:ns}) that $$\frac{1}{n}Var_{n,{\bf 0}}({\bf T}(\pi))\to \Gamma,$$ where $\Gamma$ is positive definite. Thus $Var_{n,{\bm  \theta}}({\bf T}(\pi))$ cannot be singular for $n$ large enough, and so $\mathcal{S}(\pi,.)$ is strictly convex for $n$ large enough.
\\

For verifying condition (i), note that
\begin{align}\label{eq:kl}
\notag\mathcal{S}(\pi,{\bm \theta}_1)-\mathcal{S}(\pi,{\bm \theta}_0)=&\frac{1}{n}\Big[Z_n({\bm \theta}_1)-Z_n({\bm \theta}_0)\Big]-\frac{1}{n}\Big[{\bm \theta}_1-{\bm \theta}_0\Big]^{\top}{\bf T}(\pi)\\
\stackrel{P}{\to} &Z({\bm \theta}_1)-Z({\bm \theta}_0)-({\bm \theta}_1-{\bm \theta}_0)^{\top}{\bf z}({\bm \theta})
\end{align}
where we use parts (a) and (c) of Proposition \ref{prop:sumit}. But the RHS of \eqref{eq:kl} is positive by part (f) of Proposition \ref{prop:sumit}, and so (i) holds. 
\\

Proceeding to show (ii), we have
\begin{align*}
\Big|\mathcal{S}(\pi,{\bm \theta}_1)-\mathcal{S}(\pi,{\bm \theta}_2)\Big|\le &\frac{1}{n}\Big|Z_n({\bm \theta}_1)-Z_n({\bm \theta}_2)\Big|+\frac{1}{n}|({\bm \theta}_1-{\bm \theta}_2)^{\top}{\bf T}(\pi)|\\
\le &2\|{\bm \theta}_1-{\bm \theta}_2\|_2 K\sqrt{L} ,
\end{align*}
where $K:=\max_{r\in [L]}\|f_r\|_\infty$. Thus (ii) holds with $M(\pi)=2K\sqrt{L}$, which is a constant  free of $\pi$ (and hence trivially $O_P(1)$). Thus we have verified all conditions of Lemma \ref{lem:consistb}, and so the consistency of MLE follows on noting that the unique minimizer of ${\bm \theta}\mapsto \mathcal{S}(\pi,{\bm \theta})=\frac{1}{n}[Z_n({\bm \theta})-{\bm \theta}^{\top}{\bf T}(\pi)]$ is the MLE.
\\

\item[(b)]

For proving part (b), fixing ${\bf d}=(d_1,\ldots,d_L)$ arbitrary, it suffices to show that

\begin{align*}
{\bf d}^{\top}{\bf T}(\pi)-{\bf d}^{\top}\nabla Z_n(\bm{0})=&\sum_{i=1}^n{\bf d}^{\top}{\bf f}\left(\frac{i}{n},\frac{\pi(i)}{n} \right)-\E_{\bf 0}\sum_{i=1}^n{\bf d}^{\top}{\bf f}\left(\frac{i}{n},\frac{\pi(i)}{n} \right)\\
\rightarrow &N(0,\bm{d}^{\top}\Gamma\bm{d}).
\end{align*}
 Define a function $G(.,.):[0,1]^2\mapsto \R$ by setting
\begin{align*}
G_{\bf d}(x,y):={\bf d}^{\top}{\bf f}(x,y),
\end{align*}
so that ${\bf d}^{\top}{\bf T}(\pi)=\sum_{i=1}^nG_{\bf d}\Big(\frac{i}{n},\frac{\pi(i)}{n}\Big).$
If $\bm{\theta_0}=\bm{0}$, then $\pi$ has a uniform distribution on $S_n$. 
 Applying Hoeffding's combinatorial CLT (\cite[Thm 3]{Hoeffding}), setting $$c_n(i,j):=G_{\bf d}\Big(\frac{i}{n},\frac{j}{n}\Big)-\frac{1}{n}\sum_{k=1}^n \Big[G_{\bf d}\Big(\frac{i}{n},\frac{k}{n}\Big)+G_{\bf d}\Big(\frac{k}{n},\frac{j}{n}\Big)\Big]+\frac{1}{n^2}\sum_{k,\ell=1}^n G_{\bf d}\Big(\frac{k}{n},\frac{\ell}{n}\Big)$$ we conclude that
\begin{align}\label{eq:hoeffding}
\frac{{\bf d}^{\top}{\bf T}(\pi)-{\bf d}^{\top}\nabla Z_n(\bm{0})}{\sqrt{Var_{\bf 0}({\bf d}^{\top}{\bf T}(\pi))}}\stackrel{D}{\to}N(0,1),
\end{align}
 as soon as
$$\frac{n\max_{i,j\in [n]}c_n^2(i,j)}{\sum_{i,j=1}^nc_n^2(i,j)}\to 0.$$ But the last display is immediate on noting that
\begin{align*}
\frac{n^2\max_{i,j\in [n]}c_n^2(i,j)}{\sum_{i,j=1}^nc_n^2(i,j)}\to
\frac{\sup_{x,y\in [0,1]}G_{\bf d}^2(x,y)}{\int_{[0,1]^2} G^2_{\bf d}(x,y)dx dy}<\infty.
\end{align*}
To complete the proof of part (b), using \eqref{eq:hoeffding} it suffices to show that
\begin{align}\label{eq:ns}
\frac{1}{n}Var_{n,{\bf 0}}\Big({\bf d}^{\top}{\bf T}(\pi)\Big)\to \bm{d}^{\top}\Gamma\bm{d}.
\end{align}
But this follows on invoking \cite[(Eq 10)]{Hoeffding} to note that
\begin{align*}
\frac{1}{n}Var_{n,{\bf 0}}\Big({\bf d}^{\top}{\bf T}(\pi)\Big)=&\frac{1}{n}Var_{n,{\bf 0}}\Big(\sum_{i=1}^n G_{\bf d}\Big(\frac{i}{n},\frac{\pi(i)}{n}
\Big)\Big)
=\frac{1}{n(n-1)}\sum_{i,j=1}^nc_n^2(i,j),
\end{align*}
which converges to $$ \sum_{p,q=1}^Ld_pd_q\int_{[0,1]^2}  f_p(x,y)f_q(x,y)dx dy=\bm{d}^{\top}\Gamma\bm{d},$$
as desired.

\item[(c)]

The proof of part (c) will be done by invoking Lemma \ref{lem:consist} with the following choices:
\[\mathcal{S}(\pi,{\bm \theta})=Z_n({\bm \theta})-{\bm \theta}^{\top} {\bf T}(\pi),\quad  B_1=\Gamma,\quad B_2=\Gamma.\]
In this case with ${\bf z}_n({\bm \theta}):=\nabla Z_n({\bm \theta})$ we have
\[\nabla \mathcal{S}(\pi,{\bm \theta})= \nabla Z_n({\bm \theta})-{\bf T}(\pi)={\bf z}_n({\bm \theta})-{\bf T}(\pi), \quad H(\pi,{\bm \theta})=Var_{\bm \theta}({\bf T}(\pi)),\]
and so $\mathcal{S}(\pi,.)$ is non-negative definite. Note that assumption (i) of Lemma \ref{lem:consist} holds by part (b). Also the fact that $\Gamma$ is positive definite follows immediately from the linear independence of $\{f_1,\ldots,f_r\}$. To verify condition (ii) of Lemma \ref{lem:consist} and hence complete the proof, it suffices to show that for any sequence $\{{\bf h}_n\}_{n\ge 1}\in \R^L$ converging to ${\bf h}$, we have
\begin{align}\label{eq:eee}
\frac{1}{n}Var_{n,\frac{{\bf h}_n}{\sqrt{n}}}({\bf T}(\pi))\to \Gamma.
\end{align}
%
%
%
%
To this effect, let $\{{\bf u}_n\}_{n\ge 1}$ be a sequence in $\R^L$ converging to ${\bf u}$. Then it follows by part (b) that
\[\frac{{\bf u}_n^{\top}{\bf T}(\pi)-{\bf u}_n^{\top}{\bf z}_n({\bf 0})}{\sqrt{n}}\stackrel{D}{\to} N(0,{\bf u}^{\top}\Sigma {\bf u}).\]
Also, using \cite[Prop 3.10]{Chat-thesis}
we get
 $\Big\{ \exp\Big(\frac{{\bf u}_n^{\top}{\bf T}(\pi)-{\bf u}_n^{\top}{\bf z}_n({\bf 0})}{\sqrt{n}}\Big)\Big\}_{n\ge 1}$ is uniformly integrable, and so using the last display gives
\begin{align}\label{eq:logpart1}
 \notag Z_n\Big(\frac{{\bf u}_n}{\sqrt{n}}\Big)-Z_n({\bf 0})-\frac{1}{\sqrt{n}}{{\bf u}_n}^{\top}{\bf z}_n({\bf 0})=&\log \E_{n,{\bf 0}} e^{\frac{{\bf u}_n^{\top}{\bf T}(\pi)-{\bf u}_n^{\top}{\bf z}_n({\bf 0})}{\sqrt{n}}}\\
\to& \log \E e^{ N(0, {\bf u}^{\top}\Sigma {\bf u})}=\frac{1}{2}{\bf u}^{\top}\Gamma {\bf u}.
\end{align}
In particular, choosing ${\bf u}_n={\bf u}$ for all $n$ we get
\begin{align*}
 Z_n\Big(\frac{{\bf u}}{\sqrt{n}}\Big)-Z_n({\bf 0})-\frac{1}{\sqrt{n}}{{\bf u}}^{\top}{\bf z}_n({\bf 0})\to \frac{1}{2}{\bf u}^{\top}\Gamma {\bf u}.
\end{align*}
Since the function ${\bf u}\mapsto  Z_n\Big(\frac{\bf u}{\sqrt{n}}\Big)$ is convex, and the RHS above is differentiable, the gradients converge, i.e.
\begin{align*}
\frac{{\bf z}_n\Big(\frac{{\bf u}}{\sqrt{n}}\Big)-{\bf z}_n({\bf 0})}{\sqrt{n}}\to \Gamma {\bf u}.
\end{align*}
Also, since the limiting function in the RHS above is continuous, the above convergence is uniform, i.e.~for any sequence $\{{\bf u}_n\}_{n\ge 1}$ converging to ${\bf u}$, we have
\begin{align}\label{eq:logpart2}
\frac{{\bf z}_n\Big(\frac{{\bf u}_n}{\sqrt{n}}\Big)-{\bf z}_n({\bf 0})}{\sqrt{n}}\to \Gamma {\bf u}.
\end{align}
Thus, fixing ${\bf d}\in \R^L$ and
using \eqref{eq:logpart1} with ${\bf u}_n={\bf h}_n+t{\bf d}$ and \eqref{eq:logpart2} with ${\bf u}_n={\bf h}_n$, for any $t\in \R$ we have
\begin{align*}
&\log \E_{n, {n^{-1/2}}{\bf h}_n} \exp\left[\frac{t}{\sqrt{n}}\Big({\bf d}^{\top} {\bf T}(\pi)-{\bf d}^{\top}{\bf z}_n\Big(\frac{{\bf h}_n}{\sqrt{n}}\Big)\Big)\right]\\
=&Z_n\Big(\frac{{\bf h}_n+t{\bf d}}{\sqrt{n}}\Big)-Z_n\Big(\frac{{\bf h}_n}{\sqrt{n}}\Big)-\Big(\frac{t}{\sqrt{n}}\Big){\bf d}^{\top}{\bf z}_n\Big(\frac{{\bf h}_n}{\sqrt{n}}\Big)
\to \frac{t^2}{2}{\bf d}^{\top}\Gamma {\bf d}.
\end{align*}
The above display implies that under $\P_{n,\frac{{\bf h}_n}{\sqrt{n}}}$, $$\frac{{\bf d}^{\top}{\bf T}(\pi)-{\bf d}^{\top}{\bf z}_n\Big(\frac{{\bf h}_n}{\sqrt{n}}\Big)}{\sqrt{n}}\stackrel{D}{\to} N(0, {\bf d}^{\top}\Gamma {\bf d}),$$
and the above convergence is also in moments. Consequently, we have 
%
$$\frac{1}{n}Var_{n,\frac{{\bf h}_n}{\sqrt{n}}}({\bf d}^{\top}{\bf T}(\pi))\to {\bf d}^{\top}\Gamma {\bf d}.$$
Since this holds for all ${\bf d}\in \R^L$, we conclude
$\frac{1}{n}Var_{n,\frac{{\bf h}_n}{\sqrt{n}}}({\bf T}(\pi))\to \Gamma$. This verifies \eqref{eq:eee}, and hence completes the proof of the theorem.

\end{enumerate}

\subsection{Proof of Lemma \ref{lem:CI}}

We claim that it suffices to show
\begin{align}\label{eq:suff_claim}
\widehat{A}(\hat{\bm \theta}_{PL})\stackrel{P}{\to} A({\bm \theta}_0),\quad \widehat{\Sigma}(\hat{\bm \theta}_{PL})\stackrel{P}{\to} \Sigma({\bm \theta}_0),
\end{align}
We first complete the proof of the lemma, deferring the proof of the claim. Using Theorem \ref{thm:pl} part (b) we get
$$\frac{\sqrt{n}{\bf d}^\top (\hat{\bm \theta}_{PL}-{\bm \theta}_0)}{\sqrt{ {\bf d}^\top [A({\bm \theta_0})]^{-1} \Sigma({\bm \theta}_0) [A({\bm \theta}_0)]^{-1}{\bf d}}}\stackrel{D}{\to} N\Big(0,1\Big).$$
Also, \eqref{eq:suff_claim} gives
$$ {\bf d}^\top [\widehat{A}(\hat{\bm \theta}_{PL})]^{-1} \widehat{\Sigma}(\hat{\bm \theta}_{PL}) [\widehat{A}(\hat{\bm \theta}_{PL})]^{-1}{\bf d}\stackrel{P}{\to}  {\bf d}^\top [A({\bm \theta_0})]^{-1} \Sigma({\bm \theta}_0) [A({\bm \theta}_0)]^{-1}{\bf d}.$$
Combining the above two displays along with Slutsky's theorem gives
$$\frac{\sqrt{n}{\bf d}^\top (\hat{\bm \theta}_{PL}-{\bm \theta}_0)}{  \sqrt{{\bf d}^{\top}\Big[[\widehat{A}(\hat{\bm \theta}_{PL})]^{-1}\widehat{\Sigma}(\hat{\bm \theta}_{PL})[\widehat{A}(\hat{\bm \theta}_{PL})]^{-1}\Big]{\bf d}}}\stackrel{D}{\to} N(0,1),$$
from which the construction of the confidence interval is immediate.
\\

To complete the proof, it suffices to verify \eqref{eq:suff_claim}. To this effect, setting $M:=\max_{1\le r\le L}\|f_r\|_\infty<\infty$ we have $\max_{i,j\in [n], \pi\in S_n, r\in [L]}|y_{\pi,r}(i,j)|\le 4M$. Since $\hat{\bm \theta}_{PL}\stackrel{P}{\to}{\bm \theta}_0$ (by Theorem \ref{thm:pl} part (b)), it follows from uniform continuity of the involved functions on compact domains that
$$\widehat{A}(\hat{\bm \theta}_{PL})-\widehat{A}({\bm \theta}_0)\stackrel{P}{\to}0,\quad \widehat{\Sigma}(\hat{\bm \theta}_{PL})- \widehat{\Sigma}({\bm \theta}_0)\stackrel{P}{\to}0.$$
It thus suffices to show that
$$ \widehat{A}({\bm \theta}_0)\stackrel{P}{\to} A({\bm \theta}_0),\quad \widehat{\Sigma}({\bm \theta}_0)\stackrel{P}{\to} \Sigma({\bm \theta}_0).$$
Of these, the first conclusion follows on using \eqref{eq:CI2}. For the second conclusion, use \eqref{eq:c2} to note that
$$\E_{n,{\bm \theta}_0}\widehat{\Sigma}_{pq}({\bm \theta_0})\to \Sigma_{pq}({\bm \theta_0}).$$
To complete the proof, it suffices to show that
$$Var_{n,{\bm \theta}_0}(\widehat{\Sigma}_{pq}({\bm \theta}_0))\to 0.$$
But this follows on using \eqref{eq:ejp} to note that for any finite collection of indices $\{p_1,\ldots,p_\ell\}$, the random variables $\{\pi(p_1),\ldots,\pi(p_\ell)\}$ are approximately independent.

\section{Proof of supporting lemmas}\label{sec:support}

\subsection{Proof of Proposition \ref{prop:sumit}}

\begin{enumerate}
\item[(a)] This follows on invoking \cite[Thm 1.5 (a)]{Mukherjee}  with $f={\bm \theta}^{\top}{\bf f}$ and $\theta=1$.

\item[(b)] Uniqueness of the optimization problem in part (a) follows from  \cite[Thm 1.5 (b)]{Mukherjee}. The fact that the density $\rho_{\bm \theta}$ is strictly positive everywhere follows from the form of the optimizing density in  \cite[Thm 1.5 (c)]{Mukherjee}.

\item[(c)] The convergence of the empirical measure $\nu_\pi$ follows from   \cite[Thm 1.5 (b)]{Mukherjee}. 
The second conclusion follows on noting that 
\begin{align*}
\frac{1}{n}\E_{{n,{\bm \theta}}}[{\bf T}(\pi)]=\E_{n,{\bm \theta}}[\nu_\pi({\bf f})] \stackrel{P}{\to} \mu_{\theta}({\bf f})={\bf z}({\bm \theta}),
  \end{align*}
  where we use Dominated Convergence Theorem.

\item[(d)]
 Let $\{{\bm \theta}_k\}_{k\ge 1}$ be a sequence in $\R^L$ converging to ${\bm \theta}_\infty$. We claim that 
  \begin{align}\label{eq:mu_converge}
  \mu_{{\bm \theta}_k}\stackrel{w}{\rightarrow}\mu_{{\bm \theta}_\infty}.
  \end{align} It then follows that
  $${\bf z}({\bm \theta}_k)=\mu_{{\bm \theta}_k}({\bf f})\to \mu_{{\bm \theta}_\infty}({\bf f})={\bf z}({\bm \theta}_\infty),$$
  thereby showing continuity.
  Proceeding to show \eqref{eq:mu_converge}, let $\mu_\infty$ be any limit point of the sequence of measures $\{\mu_{{\bm \theta}_k}\}_{k\ge 1}$, which exists by tightness of $\mathcal{M}$. Using continuity of $Z(.)$ and lower semi continuity of $D(.||u)$, we have
  \begin{align*}
  Z({\bm \theta}_\infty)=\lim_{k\to\infty}Z({\bm \theta}_k)=&\lim_{k\to\infty}\Big \{ {\bm \theta}_k \mu_{{\bm \theta}_k}({\bf f})-D(\mu_{{\bm \theta}_k}||u)\Big\}
  \le {\bm \theta}_\infty \mu({\bf f})-D(\mu_\infty||u).
  \end{align*}
  But from parts (a) and (b) we have
  $$Z({\bm \theta}_\infty)=\sup_{\mu\in \mathcal{M}}\Big\{{\bm \theta}_\infty \mu({\bf f})-D(\mu||u)\Big\}={\bm \theta}_\infty \mu_{{\bm \theta}_\infty}({\bf f})-D(\mu_{{\bm \theta}_\infty}||u).$$
Comparing last two displays along with the uniqueness of optimizer from part (b) gives $\mu_{{\bm \theta}_\infty}=\mu_\infty$. Thus we have shown \eqref{eq:mu_converge}, and so continuity of ${\bf z}(.)$ follows.
\\

\item[(e)]
Using part (c) we get
  \begin{align*}
  \frac{\nabla Z_n({\bm \theta})}{n}=\frac{1}{n}\E_{{n,{\bm \theta}}}[{\bf T}(\pi)]\stackrel{P}{\to} {\bf z}({\bm \theta}).
  \end{align*}
   Since the function $Z(.)$ is limit of convex functions, it is convex, and hence differentiable a.s.~Using part (a) and the display above, we then have
  $$\nabla Z({\bm \theta})\stackrel{a.s.}{= }{\bf z}({\bm \theta}).$$ Finally, since ${\bf z}(.)$ is continuous (by part (d) above), it follows that $Z(.)$ is differentiable everywhere, and the above equality holds everywhere.\\
  
  \item[(f)]

To show strict convexity, assume by way of contradiction that there exists $\theta_0\ne \theta_1$ such that 
\begin{align}\label{eq:fail}
Z({\bm \theta}_1)=Z({\bm \theta}_0)+({\bm \theta}_1-{\bm \theta}_0)^{\top} {\bf z}({\bm \theta_0}).
\end{align}
Define a function $\phi:[0,1]\mapsto \R$ by setting
$$\phi(t):=Z\Big((1-t){\bm \theta}_0+t{\bm \theta}_1\Big)=Z\Big({\bm \theta}_0+t({\bm \theta}_1-{\bm \theta}_0)\Big).$$
Then $\phi$ is convex, and 
$$\phi(0)=Z({\bm \theta}_0), \quad \phi(1)=Z({\bm \theta}_1), \quad \phi'(t)=({\bm \theta}_1-{\bm \theta}_0)^{\top}{\bf z}\Big({\bm \theta}_0+t({\bm \theta}_1-{\bm \theta}_0)\Big).$$
Consequently,  \eqref{eq:fail} translates into $\phi(1)=\phi(0)+\phi'(0)$. Convexity of $\phi$ then implies that $\phi(t)=\phi(0)+t\phi'(0)$ for all $t\in [0,1]$. Differentiating, we see that $\phi'(t)=\phi'(0)$ for all $t\in [0,1]$, and so 
\begin{align}\label{eq:der=0}
({\bm \theta}_1-{\bm \theta}_0)^{\top}{\bf z}({\bm \theta}_1)=\phi'(1)=\phi'(0)=({\bm \theta}_1-{\bm \theta}_0)^{\top}{\bf z}({\bm \theta}_0).
\end{align}
Now, optimality of $\mu_{{\bm \theta}_1}$ gives
$${\bm \theta}_1^{\top} \mu_{{\bm \theta}_1}({\bf f})-D(\mu_{{\bm \theta}_1}||u)\ge {\bm \theta}_1^{\top}\mu_{{\bm \theta}_0}({\bf f})-D(\mu_{{\bm \theta}_0}||u),$$
which implies
\begin{align}\label{eq:opt1}
{\bm \theta}_1^{\top}\Big({\bf z}({\bm \theta}_1)-{\bf z}({\bm \theta}_0)\Big)\ge D(\mu_{{\bm \theta}_1}||u)-D(\mu_{{\bm \theta}_0}||u).
\end{align}
Similarly, optimality of $\mu_{{\bm \theta}_0}$ gives
\begin{align}\label{eq:opt2}
{\bm \theta}_0^{\top}\Big({\bf z}({\bm \theta}_1)-{\bf z}({\bm \theta}_0)\Big)\le D(\mu_{{\bm \theta}_1}||u)-D(\mu_{{\bm \theta}_0}||u).
\end{align}
Combining \eqref{eq:opt1} and \eqref{eq:opt2} gives
\begin{align*}
{\bm \theta}_1^{\top}\Big({\bf z}({\bm \theta}_1)-{\bf z}({\bm \theta}_0)\Big)\ge D(\mu_{{\bm \theta}_1}||u)-D(\mu_{{\bm \theta}_0}||u)\ge {\bm \theta}_0^{\top}\Big({\bf z}({\bm \theta}_1)-{\bf z}({\bm \theta}_0)\Big).
\end{align*}
But by \eqref{eq:der=0} the LHS and RHS of the above display are the same, and so
\begin{align*}
{\bm \theta}_1^{\top}\Big({\bf z}({\bm \theta}_1)-{\bf z}({\bm \theta}_0)\Big)= D(\mu_{{\bm \theta}_1}||u)-D(\mu_{{\bm \theta}_0}||u)= {\bm \theta}_0^{\top}\Big({\bf z}({\bm \theta}_1)-{\bf z}({\bm \theta}_0)\Big).
\end{align*}
This implies that equality must hold in both \eqref{eq:opt1} and \eqref{eq:opt2}, and so uniqueness of optimizer (from part (b)) implies  $\mu_{{\bm \theta}_1}=\mu_{{\bm \theta}_0}$. This, along with the form of the density in  \cite[Thm 1.5 (c)]{Mukherjee} gives
\begin{align*}
\rho_{{\bm \theta}_1}(x,y)=&\exp\Big({\bm \theta}_1^{\top}{\bf f}(x,y)+a_{{\bm \theta}_1}(x)+b_{{\bm \theta}_1}(y)\Big)\\
=&\exp\Big({\bm \theta}_0^{\top}{\bf f}(x,y)+a_{{\bm \theta}_0}(x)+b_{{\bm \theta}_0}(y)\Big)=\rho_{{\bm \theta}_0}(x,y),
\end{align*}
for some functions $a_{\bm \theta}, b_{\bm \theta}$ from $[0,1]$ to $\R$.
Taking $\log$ and simplifying gives
\begin{align}\label{eq:sumitc}
({\bm \theta}_1-{\bm \theta}_0)^{\top}{\bf f}(x,y)=a_{{\bm \theta}_0}(x)+ b_{{\bm \theta}_0}(y)-a_{{\bm \theta}_1}(x)-b_{{\bm \theta}_1}(y).
\end{align}
Integrating over $y\in [0,1]$ and using the fact that $f_r\in \mathcal{C}$ for all $r\in [L]$ we get
$$a_{{\bm \theta}_1}(x)-a_{{\bm \theta}_0}(x)=\int_0^1 [b_{{\bm \theta}_0}(y)-b_{{\bm \theta}_1}(y)]dy.$$
Thus the function $ a_{{\bm \theta}_1}(.)-a_{{\bm \theta}_0}(.)$ is constant a.s.~A similar calculation by integrating over $x$ gives that $b_{{\bm \theta}_1}(.)-b_{{\bm \theta}_0}(.)$ is constant a.s., and so using \eqref{eq:sumitc} it follows that $({\bm \theta}_1-{\bm \theta}_0)^{\top}{\bf f}(x,y)$ is a constant a.s.~But the assumption $f_r\in \mathcal{C}$ then forces this constant to be $0$, which is a contradiction to the assumption that ${\bf f}$ is linearly independent. Thus \eqref{eq:fail} cannot hold, and so $Z(.)$ is strictly convex, as claimed.

\end{enumerate}

%


\subsection{Proof of Lemma \ref{lem:conditional_mean}}
Using the conditional distribution of  $(\pi(i),\pi(j))$ given $\{\pi(\ell), \ell\ne i,j\}$ as specified in \eqref{eq:pi_condition}, we get
\begin{align*}
&\E_{n,{\bm \theta}}\Big({\bf C}_{ij}(\pi)\Big| \pi(\ell)=\sigma(\ell), \ell \ne i,j\Big)\\
=&\frac{e^{{\bm \theta}^{\top} {\bf y}_{\sigma}(i,j)}}{1+e^{{\bm \theta}^{\top} {\bf y}_{\sigma}(i,j)}}\frac{1}{1+e^{{\bm \theta}^{\top}{\bf y}_{\sigma}(i,j)} }{\bf y}_{\sigma}(i,j)+\frac{1}{1+e^{{\bm \theta}^{\top} {\bf y}_{\sigma}(i,j)}}\frac{1}{1+e^{-{\bm \theta}^\top {\bf y}_{\sigma}(i,j)}}[-{\bf y}_{\sigma}(i,j)] \\
=&\frac{e^{{\bm \theta}^\top {\bf y}_{\sigma}(i,j)}}{\Big(1+e^{{\bm \theta}^\top {\bf y}_{\sigma}(i,j)}\Big)^2}\Big[{\bf y}_\sigma(i,j)-{\bf y}_\sigma(i,j)\Big]={\bf 0}.
\end{align*}
The desired conclusion follows on noting that $C_{ij}(\pi,{\bf d})={\bf d}^\top {\bf C}_{ij}(\pi)$.

\subsection{Proof of Lemma \ref{lem:moment}}

To begin, a direct expansion gives
\begin{align}\label{eq:direct}
&\Big|\E_{n,{\bm \theta_0}} \Big(\sum_{(i,j)\in \mathcal{E}_n}C_{ij}(\pi,{\bf d})\Big)^k-\E D_n^k\Big|\nonumber\\&\le  \sum_{(i_1,j_1),\ldots,(i_k,j_k)\in \mathcal{E}_n}\left| \E_{n,{\bm \theta_0}}\Big[ \prod_{a=1}^k C_{i_aj_a}(\pi)\Big] -\E \Big[\prod_{a=1}^k D_{i_aj_a}\Big]\right|.
\end{align}
Given a collection of $k$ pairs $(i_1,j_1),\ldots, (i_k,j_k)$ in $\mathcal{E}_n$, let $H({\bf i},{\bf j})$ denote the labelled multigraph formed by the edges $(i_1,j_1),\ldots,(i_k,j_k)$. The graph $H({\bf i},{\bf j})$ contains multiple edges, but no loops. Depending on the graph $H({\bf i},{\bf j})$, we split the argument into the following two cases:

\begin{itemize}
\item{\bf Case 1:} $H=H({\bf i},{\bf j})$ has an isolated edge, i.e. it has an edge $(i_1,j_1)$ which does not intersect any of the other edges $\{(i_a,j_a), 2\le a\le k\}$.
\\

In this case we will show that
\begin{align}\label{eq:case1}
 \E_{n,{\bm \theta_0}} \Big[\prod_{a=1}^k C_{i_aj_a}(\pi)\Big] =0,\quad \E \Big[\prod_{a=1}^k D_{i_aj_a}\Big]=0.
 \end{align}
Indeed, invoking Lemma \ref{lem:conditional_mean} gives
\begin{align*}
\E_{n,{\bm \theta_0}}\Big[ \prod_{a=1}^k C_{i_aj_a}(\pi)\Big]= \E_{n,{\bm \theta_0}}\Big[\E_{n,{\bm \theta_0}} \Big( C_{i_1j_1}(\pi)\Big| \pi(\ell), \ell\ne i_1,j_1\Big)  \prod_{a=2}^k C_{i_aj_a}(\pi)\Big]=0,
\end{align*}
thus verifying the first equality in \eqref{eq:case1}. Proceeding to show \eqref{eq:case1}, again using Lemma \ref{lem:conditional_mean} along with the construction of $\{D_{ij},(i,j)\in \mathcal{E}_n\}$ gives
$$\E \Big[D_{i_1j_1} D_{\ell r}\Big]=\E_{n,{\bm \theta_0}}\Big[ C_{i_1j_1}(\pi) C_{\ell r}(\pi)\Big]=0\text{ whenever }(i_1,j_1)\cap (\ell,r)=\phi.$$ Since $\{D_{ij}, (i,j)\in \mathcal{E}_n\}$ is multivariate Gaussian, we conclude that $D_{i_1j_1}$ is independent of $\{D_{\ell r}, (\ell,r)\cap (i_1,j_1)=\phi\}$. Consequently, we have
 $$ \E\Big(D_{i_1j_1}|D_{\ell r}, (\ell,r)\cap (i_1,j_1)=\phi\Big)= \E D_{i_1j_1}=0,$$
which gives 
$$\E\Big[ \prod_{a=1}^k D_{i_aj_a}\Big]= \E\Big[\E \Big(D_{i_1 j_1}\Big|D_{\ell r}, (\ell, r)\cap (i_1,j_1)=\phi\Big) \prod_{a=2}^k D_{i_a j_a}\Big]=0.$$
This verifies the second equality of \eqref{eq:case1}, and hence completes the proof of case 1.
\\

\item{\bf Case 2:} $H=H({\bf i},{\bf j})$ has more than $\frac{3k}{2}$ vertices.
\\

In this case we will show that there is an isolated edge. Consequently, using case 1 above we will get
$$ \E_{n,{\bm \theta_0}} \Big[\prod_{a=1}^k C_{i_aj_a}(\pi)\Big] =0,\quad \E \Big[\prod_{a=1}^k D_{i_aj_a}\Big]=0.$$

Proceeding to find an isolated edge,  let $\cup_{b=1}^sH_b$ denote the decomposition of the graph $H$ into $s$ connected components. By way of contradiction, assume that no isolated edge exists. Then each component $H_b$ has $|V(H_b)|\ge 3, |E(H_b)|\ge 2$, which gives
\begin{align}\label{eq:graph1}
k=|E(H)|=\sum_{b=1}^s |E(H_b)|\ge 2s.
\end{align}
Also since each $H_b$ is connected, we have $|E(H_b)|\ge |V(H_b)|-1$. On adding over $b$ we get
\begin{align}\label{eq:graph2}
k=|E(H)|\ge |V(H)|-s>\frac{3k}{2}-s.
\end{align}
Combining \eqref{eq:graph1} and \eqref{eq:graph2} we have $s\le \frac{k}{2}$ and $s>\frac{k}{2}$, a contradiction. This completes the proof of case 2.
\\

\item{\bf Case 3:} $H({\bf i},{\bf j})$ has exactly $\frac{3k}{2}$ vertices, and no edge in $H$ is isolated.
\\

In this case we will show that $H$ is a disjoint union of two stars, i.e. each component of $H$ is a connected graph with two edges and three vertices. Consequently, invoking \eqref{eq:ejp} we get
\begin{align*}
&\E_{n,{\bm \theta_0}}\Big[ \prod_{a=1}^k C_{i_aj_a}(\pi)\Big]\\
 =&\prod_{b=1}^{\frac{k}{2}} \E_{n,{\bm \theta_0}} \Big[\prod_{(i_a,j_a)\in E(H_b)} C_{i_a j_a}(\pi)\Big]+o(1)\\
=&\prod_{b=1}^{\frac{k}{2}} \E \Big[\prod_{(i_a,j_a)\in E(H_b)} D_{i_a j_a}\Big]+o(1)\text{ [Since $H_b$ has exactly two edges }]\\
=&\E \Big[\prod_{a=1}^k D_{i_aj_a}\Big]+o(1)\text{ [By independence of uncorrelated Gaussians}].
\end{align*}

Proceeding to verify that $H$ is indeed a disjoint union of two stars, first note that since $H$ has $\frac{3k}{2}$ vertices, $k$ must be even. Also, the same proof as in case 2 shows that the number of connected components of $H$ must be exactly $s=\frac{k}{2}$. Also the assumption that $H$ has no isolated edges gives $|V(H_b)|\ge 3, |E(H_b)|\ge 2$. This gives
$$\frac{3k}{2}=|V(H)|=\sum_{b=1}^{\frac{k}{2}}|V(H_b)|\ge \frac{3k}{2},\quad k=|E(H)|=\sum_{b=1}^{\frac{k}{2}}|E(H_b)|\ge k,$$
and so each component of $H$ has exactly $3$ vertices and $2$ edges, as desired. This completes the proof of case 3.

\end{itemize}

Combining the three cases above, and noting that the number of ways to choose ${\bf i}, {\bf j}$ such that $H({\bf i},{\bf j})$ has $r$ vertices is $O(n^r)$, \eqref{eq:direct} gives

\begin{align*}
\Big|\E_{n,{\bm \theta_0}} \Big(\sum_{(i,j)\in \mathcal{E}_n}C_{ij}(\pi,{\bf d})\Big)^k-\E D_n^k\Big|\le  \sum_{1\le r<\frac{3k}{2}} O(n^r) +o(n^{\frac{3k}{2}})1\{k\text{ is even}\},
\end{align*}
from which the conclusion of Lemma \ref{lem:moment} follows.

\subsection{Proof of Lemma \ref{lem:consist}} 


To begin, note that the function $\mathcal{S}(\pi,.)$ is convex.
Throughout the proof, for any $\delta>0$ and ${\bm \theta}_1\in \R^L$ set $$B({\bm \theta}_1,\delta):=\{{\bm \theta}\in \R^L:\|{\bm \theta}-{\bm \theta}_1\|_2\le \delta\}, \quad \partial B({\bm \theta}_1,\delta):=\{{\bm \theta}\in \R^L:\|{\bm \theta}-{\bm \theta}_1\|_2= \delta\}.$$
We now break the proof into the following steps:

\begin{itemize}
\item 
For every $K>0$ we have
$$\lim_{n\to\infty}\P_{n,{\bm \theta}_0}\left(\text{ There exists a unique minimizer in } B\Big({\bm \theta}_0,\frac{K}{\sqrt{n}}\Big)\right)=1.$$

\begin{proof}
Fixing $K>0$, the function $ \mathcal{S}(\pi,.)$ is continuous, and hence achieves its minimum on the compact set $B\Big({\bm \theta}_0,\frac{K}{\sqrt{n}}\Big)$. We now argue uniqueness of minimizer. Indeed, if $\theta_{n,K,1}$ and $\theta_{n,K,2}$ are two minimizers of $\mathcal{S}(\pi,.)$ in $ B\Big({\bm \theta}_0,\frac{K}{\sqrt{n}}\Big)$, then using the convexity of $\mathcal{S}(\pi,.)$ it must be constant on the line joining $\theta_{n,K,1}$ and $\theta_{n,K,2}$, and hence the gradient must vanish on that line. Thus there exists a random variable ${\bm \xi}_{n,K}\in B\Big({\bm \theta}_0, \frac{K}{\sqrt{n}}\Big)$ such that the matrix $H(\pi,{\bm \xi}_{n,K})$ is singular. But then using assumption (ii) it follows that $n^{-1}H(\pi,{\bm \xi}_{n,K})\stackrel{P}{\to}B_2$ with $B_2$ non singular, which is a contradiction. Thus with probability tending to $1$, there is a unique minimizer $\theta_{n,K}$ on the set $ B\Big({\bm \theta}_0,\frac{K}{\sqrt{n}}\Big)$.
\\
\end{proof}


\item
 $$\limsup\limits_{K\to\infty}\limsup\limits_{n\to\infty}\P_{n,{ \bm \theta}_0}\Big( \mathcal{S}(\pi,{\bm \theta}_0)\ge  \inf_{{\bm \theta}\in  \partial B\Big({\bm \theta}_0,\frac{K}{\sqrt{n}}\Big)}\mathcal{S}(\pi,{\bm \theta})\Big)=0.$$

\begin{proof}

Let ${\bm \theta}_{n,K}$ be a global minimizer of $\mathcal{S}(\pi,.)$ over the compact set $ \partial B\Big({\bm \theta}_0,\frac{K}{\sqrt{n}}\Big)$.
Setting ${\bf u}_{n,K}:=\frac{\sqrt{n}}{K}({\bm \theta}_{n,K}-{\bm \theta}_0)$ we have $\|{\bf u}_{n,K}\|=1$. A second order Taylor's expansion gives
\begin{align*}
\notag \mathcal{S}(\pi,{\bm \theta}_{n,K})-\mathcal{S}(\pi,{\bm \theta}_0)=&({\bm \theta}_{n,K}-{\bm \theta}_0)^\top\nabla \mathcal{S}(\pi,{\bm \theta}_0)+\frac{1}{2}({\bm \theta}_{n,K}-{\bm \theta}_0)^{\top}H(\pi,{\bm \xi}_{n,K})({\bm \theta}_{n,K}-{\bm \theta}_0)\\
\notag=&\frac{K}{\sqrt{n}}{\bf u}_{n,K}^{\top} \nabla \mathcal{S}(\pi,{\bm \theta}_0)+\frac{K^2}{2n} {\bf u}_{n,K}^{\top}H(\pi,{\bm \xi}_{n,K}){\bf u}_{n,K}\\
\ge &-K\frac{ \|\nabla \mathcal{S}(\pi,{\bm \theta}_0)\|_2}{\sqrt{n}}+\frac{K^2}{2n}\lambda_1(H(\pi,{\bm \xi}_{n,K})).
\end{align*}
This gives
\begin{align}\label{eq:contradict2}
\notag&\P_{n,{\bm \theta}_0}\Big(\mathcal{S}(\pi,{\bm \theta}_0)\ge \inf_{{\bm \theta}\in \partial B\Big({\bm \theta}_0,\frac{K}{\sqrt{n}}\Big)}\mathcal{S}(\pi,{\bm \theta})\Big)\\
\notag=&\P_{n,{\bm \theta}_0}\Big(\mathcal{S}(\pi,{\bm \theta}_0)\ge \mathcal{S}(\pi,{\bm \theta}_{n,K})\Big)\\
\le &\P_{n,{\bm \theta}_0}\left(\frac{ \|\nabla \mathcal{S}(\pi,{\bm \theta}_0)\|_2}{\sqrt{n}}\ge \frac{K}{2n}\lambda_1\Big(H(\pi,{\bm \xi}_{n,K})\Big)\right),
\end{align}
where $\lambda_1$ denotes the minimum eigenvalue of a matrix. Using (ii), we have
\[n^{-1}H(\pi, {\bm \xi}_{n,K})\stackrel{P}{\to}B_2 \text{             }\Rightarrow\text{            } n^{-1}\lambda_1(H(\pi,{\bm \xi}_{n,K}))\stackrel{P}{\to}\lambda_1(B_2)>0.\]
On the other hand, using (i) we have $\|\nabla \mathcal{S}(\pi,{\bm \theta}_0)\|_2=O_P(\sqrt{n})$. The desired conclusion follows from these two observations, along with \eqref{eq:contradict2}.
\\

\end{proof}

\item
$\lim\limits_{n\to\infty}\P_{n,{\bm \theta}_0}(\text{There exists a unique global minimizer $\hat{\bm \theta}_n$ (say) in }\R^L)=1.$ 

Further, $\sqrt{n}(\hat{\bm \theta}-{\bm \theta}_0)=O_P(1)$.

\begin{proof}
Using the previous step, there exists a sequence of positive reals $\{K_n\}_{n\ge 1}$ diverging to $\infty$, such that the following hold:
\begin{align*}
\lim_{n\to\infty}&\P_{n,{\bm \theta}_0}\left(\text{There exists a unique minimizer in }B\Big({\bm \theta}_0,\frac{K_n}{\sqrt{n}}\Big)\right)=1,\\
\lim_{n\to\infty}&\P_{n,{ \bm \theta}_0}\Big( \mathcal{S}(\pi,{\bm \theta}_0)< \inf_{{\bm \theta}\in \partial B\Big({\bm \theta}_0,\frac{K_n}{\sqrt{n}}\Big)}\mathcal{S}(\pi,{\bm \theta})\Big)=1.
\end{align*}
Without loss of generality we work on the intersection of the two sets above. 
 Let $\hat{\bm \theta}_{n}$ be the unique global minimizer of $\mathcal{S}(\pi,.)$ over $B\Big({\bm \theta}_0,\frac{K_n}{\sqrt{n}}\Big)$. Then $\hat{\bm \theta}_{n}$ is an interior point of $B\Big({\bm \theta}_0,\frac{K_n}{\sqrt{n}}\Big)$, and so $\nabla \mathcal{S}(\pi,\hat{\bm \theta}_{n})={\bf 0}$. 
 Since $\mathcal{S}(\pi,.)$ is convex, $\hat{\bm \theta}_n$ is a global minimizer.
It remains to show that there are no other global minimizers. But this follows on using convexity to note that existence of other global minimizers will force existence of other global minimizers in the set $B\Big({\bm \theta}_0, \frac{K}{\sqrt{n}}\Big)$, which is a contradiction.
 \\
 

To verify the second statement, let $\{\alpha_n\}_{n\ge 1}$ be a sequence of positive reals diverging to $\infty$. Then note that the above proof goes through verbatim by replacing $K_n$ by $K_n':=\min(K_n,\sqrt{\alpha_n})$, and so repeating the above prove gives $\sqrt{n}\|\hat{\bm \theta}-{\bm \theta}_0\|<K_n'=o(\alpha_n)$. Since this holds for an arbitrary sequence of positive reals $\{\alpha_n\}_{n\ge 1}$ diverging to $\infty$, it follows that $\sqrt{n}(\hat{\bm \theta}_n-{\bm \theta}_0)=O_P(1)$, as desired.

\end{proof}

\item 
$\hat{\bm \theta}_n$ is the unique solution of $\nabla \mathcal{S}(\pi,{\bm \theta})={\bf 0}$.

\begin{proof}
Since $\hat{\bm \theta}_n$ is the global minimizer of a twice differentiable function, it follows that  $\nabla \mathcal{S}(\pi,\hat{\bm \theta}_n)={\bf 0}$. Also, uniqueness follows on noting that $\mathcal{S}(\pi,.)$ is convex, and $\hat{\bm \theta}_n$ is the unique global minimizer.
\end{proof}

\item
$\sqrt{n}(\hat{\bm \theta}_n-{\bm \theta}_0)\stackrel{D}{\to}N({\bf 0},B_2^{-1}B_1 B_2^{-1})$.

\begin{proof}
A Taylor's series expansion gives
\begin{align*}
&{\bf 0}=\nabla \mathcal{S}(\pi,\hat{\bm \theta}_n)
=\nabla \mathcal{S}(\pi,{\bm \theta}_0)+\int_0^1 H\Big(\pi, {\bm \theta}_0+t(\hat{\bm \theta}_n-{\bm \theta}_0)\Big) (\hat{\bm \theta}_n-{\bm \theta}_0)dt,
\end{align*}
which gives
\begin{align*}
-\frac{1}{\sqrt{n}} \nabla \mathcal{S}(\pi,{\bm \theta}_0)=\int_0^1 \left[\frac{1}{n} H\Big(\pi, {\bm \theta}_0+t(\hat{\bm \theta}_n-{\bm \theta}_0)\Big)\right] \sqrt{n}(\hat{\bm \theta}_n-{\bm \theta}_0)dt.
\end{align*}
The LHS above converges in distribution to $N({\bf 0}, B_1)$ using (i). For the RHS,
using (ii) along with the fact that $\sqrt{n}(\hat{\bm \theta}_n-{\bm \theta}_0)=O_P(1)$, for every $t\in [0,1]$ we get \[\frac{1}{n}H\Big(\pi, {\bm \theta}_0+t(\hat{\bm \theta}_n-{\bm \theta}_0)\Big) \stackrel{P}{\to}B_2.\] Using DCT, letting $n\to\infty$ gives $\sqrt{n}(\hat{\bm \theta}_n-{\bm \theta}_0)\stackrel{D}{\to}{\bf W}$, where
$B_2 {\bf W}\sim N({\bf 0}, B_1)$, or equivalently, ${\bf W}\sim N({\bf 0},B_2^{-1}B_1B_2^{-1})$. This completes the proof of the lemma.

\end{proof}

\end{itemize}

\subsection{Proof of Lemma \ref{lem:consistb}}


As in the proof of Lemma \ref{lem:consist}, for every $\delta>0$, set $$B({\bm \theta}_0,\delta):=\{{\bm \theta}\in \R^L:\|{\bm \theta}-{\bm \theta}_0\|_2\le \delta\},\quad \partial B({\bm \theta}_0,\delta):=\{{\bm \theta}\in \R^L:\|{\bm \theta}-{\bm \theta}_0\|_2= \delta\}.$$ We break the proof of the lemma into several steps.

\begin{itemize}
\item
For every ${\bm \theta}_1\ne {\bm \theta}_0$ there exists $\eta:=\eta_{{\bm \theta}_1}>0$ such that
$$\P_{n,{\bm \theta}_0}\Big(\mathcal{S}(\pi,{\bm \theta}_0)<\inf_{{\bm \theta}\in  B({\bm \theta}_1,\eta) }\mathcal{S}(\pi,{\bm \theta})\Big)\to 1.$$

\begin{proof}
Assume by way of contradiction that no such $\eta>0$ exists. Then there exists a sequence of positive reals $\{\eta_n\}_{n\ge 1}$, such that for along a subsequence in $n$ we have
\begin{align}\label{eq:contra1}
\limsup_{n\to\infty}\P_{n,{\bm \theta}_0}\Big(\mathcal{S}(\pi,{\bm \theta}_0)\ge \inf_{{\bm \theta}\in B({\bm \theta}_1,\eta_n) }\mathcal{S}(\pi,{\bm \theta})\Big)>0.
\end{align}
Using assumption (ii), we get
$$\Big|\mathcal{S}(\pi,{\bm \theta}_1)-\inf_{{\bm \theta}\in B({\bm \theta}_1,\eta_n)}\mathcal{S}(\pi,{\bm \theta})\Big|\le M(\pi) \eta_n.$$
Also, using assumption (i), without loss of generality we can work on the set
$$\mathcal{S}(\pi,{\bm \theta}_0)+\varepsilon<\mathcal{S}(\pi,{\bm \theta}_1).$$
On the intersection of the two sets above, we have 
$$ \mathcal{S}(\pi,{\bm \theta}_0)\ge \mathcal{S}(\pi,{\bm \theta}_1)-M(\pi)\eta_n>\mathcal{S}(\pi,{\bm \theta}_0)+\varepsilon-M(\pi)\eta_n\Rightarrow M(\pi)>\frac{\varepsilon}{\eta_n}.$$ This gives
\begin{align*}
 \P_{n,{\bm \theta}_0}\Big(\mathcal{S}(\pi,{\bm \theta}_0)\ge \inf_{{\bm \theta}\in B({\bm \theta}_1,\eta_n) }\mathcal{S}(\pi,{\bm \theta})\Big)
\le o(1)+\P_{n,{\bm \theta}_0}\Big(M(\pi)>\frac{\varepsilon}{\eta_n}\Big)=o(1),
\end{align*}
where the last equality uses the fact that $M(\pi)=O_P(1)$. But this is a contradiction to \eqref{eq:contra1}, and so the proof is complete.
\end{proof}

\item
$$\P_{n,{\bm \theta}_0}(E_{n,\delta})\to 1,\text{ where }E_{n,\delta}:=\Big\{\pi\in S_n:\mathcal{S}(\pi,{\bm \theta}_0)<\inf_{{\bm \theta}\in \partial B({\bm \theta}_0,\delta)}\mathcal{S}(\pi,{\bm \theta})\Big\}.$$

\begin{proof}
For every ${\bm \theta}_1\ne {\bm \theta}_0$ let $\eta_{{\bm \theta}_1}$ be as constructed in the step above. The the collection of balls
$$\left\{B\Big({\bm \theta}_1, \eta_{{\bm \theta}_1}\Big), {\bm \theta}_1\in \partial B({\bm \theta}_0,\delta)\right\}$$
is a covering of the compact set $\partial B({\bm \theta}_0,\delta)$, and so there exists a positive integer $N<\infty$ and points $\{{\bm \theta}_a\}_{1\le a\le N}\in \partial B({\bm \theta}_0,\delta)$ such that
$$\partial B({\bm \theta}_0,\delta)\subseteq \cup_{a=1}^N B\Big({\bm \theta}_a, \eta_{{\bm \theta}_a}\Big).$$
Consequently, we have
\begin{align*}
\inf_{{\bm \theta}\in \partial B({\bm \theta}_0,\delta)}\mathcal{S}(\pi,{\bm \theta})\ge \min_{1\le a\le N} \inf_{{\bm \theta}\in B\Big({\bm \theta}_a, \eta_{{\bm \theta}_a}\Big)} \mathcal{S}(\pi, {\bm \theta}).
\end{align*}
But by the step above (and choice of $\eta_{{\bm \theta}_a}$), we have
$$\P_{n,{\bm \theta}_0}\Big(\mathcal{S}(\pi,{\bm \theta}_0)<\min_{1\le a\le N} \inf_{{\bm \theta}\in B\Big({\bm \theta}_a, \eta_{{\bm \theta}_a}\Big)} \mathcal{S}(\pi, {\bm \theta})\Big)\to 1.$$
Combining the above two displays we get
$$\P_{n,{\bm \theta}_0}\Big(\mathcal{S}(\pi,{\bm \theta}_0)<\inf_{{\bm \theta}\in \partial B({\bm \theta}_0,\delta)}\mathcal{S}(\pi,{\bm \theta})\Big)\to 1,$$
which is what we want.
\end{proof}

\item
 $$\P_{n,{\bm \theta}_0}(\widetilde{E}_{n,\delta})\to 1,\text{ where }\widetilde{E}_{n,\delta}:=\Big\{\pi\in S_n: \nabla \mathcal{S}(\pi, {\bm \theta})={\bf 0} \text{ has a root in }B({\bm \theta}_0,\delta)\Big\}.$$

\begin{proof}
To begin, note that global minimizers of the function $\mathcal{S}(\pi,.)$ over the set $B({\bm \theta}_0,\delta)$ exist, because the map $\mathcal{S}(\pi,.)$ is continuous by assumption (ii), and the set $B({\bm \theta}_0,\delta)$ is compact. Using the above step, without loss of generality we can work on the set $E_{n,\delta}$. On this set, any global minimizer ${\bm \theta}_{n,\delta}$ of $\mathcal{S}(\pi,.)$  over the set $B({\bm \theta}_0,\delta)$ must be attained in the interior of $B({\bm \theta}_0,\delta)$, and hence must satisfy $\nabla \mathcal{S}(\pi, {\bm \theta}_{n,\delta})={\bf 0}$. Thus we have shown that $E_{n,\delta}\subseteq \widetilde{E}_{n,\delta}$, and so we are done.
\end{proof}

\item
Under $\P_{n,{\bm \theta}_0}$ we have 
$\hat{\bm \theta}_n\stackrel{P}{\to}{\bm \theta}_0$.

\begin{proof}
Using the above step, without loss of generality we can work on the set $\widetilde{E}_{n,1}$. Thus we have the existence of  ${\bm \theta}_{n,1}\in B({\bm \theta}_0,1)$ such that $\nabla \mathcal{S}(\pi, {\bm \theta}_{n,1})={\bf 0}$. Since the function $\mathcal{S}(\pi,.)$ is strictly convex, it follows that ${\bm \theta}_{n,1}$ is the unique global minimizer of  the function  $\mathcal{S}(\pi,.)$ over $\R^L$. To complete the proof, it suffices to show that for any $\delta>0$ we have 
$$\P_{n,{\bm \theta}_0}(\|{\bm \theta}_{n,1}-{\bm \theta}_0\|_2>\delta)\to 0.$$
Again using the above step,  without loss of generality we can also work on the set $\widetilde{E}_{n,\delta}$ (note that $\widetilde{E}_{n,\delta}\subseteq \widetilde{E}_{n,1}$). On this set we have ${\bm \theta}_{n,1}={\bm \theta}_{n,\delta}$, and so
$$\|{\bm \theta}_{n,1}-{\bm \theta}_0\|_2=\|{\bm \theta}_{n,\delta}-{\bm \theta}_0\|_2\le \delta.$$
The desired consistency is immediate.
\end{proof}

\end{itemize}

\begin{acks}[Acknowledgments]
The first author was supported by NSF Grants DMS-1712037 and DMS-2113414.

\end{acks}

\bibliographystyle{imsart-number} 
\bibliography{bibliography.bib}       

%
%
%
%
%
%
%
%
%
%
%
%
%
%

\end{document}